\documentclass[11pt]{elsarticle}
\usepackage{amsmath,amssymb,amsthm}
\usepackage{hyperref}

\theoremstyle{plain}
\newtheorem{theorem}{Theorem}[section]
\newtheorem{lemma}[theorem]{Lemma}

\newtheorem{proposition}[theorem]{Proposition}
\newtheorem{corollary}[theorem]{Corollary}

\theoremstyle{definition}
\newtheorem{definition}[theorem]{Definition}
\newtheorem{remark}[theorem]{Remark}
\newtheorem{example}[theorem]{Example}

\newcommand{\dda}{\mathord{\mbox{\makebox[0pt][l]{\raisebox{-.4ex}{$\downarrow$}}$\downarrow$}}}
\newcommand{\dua}{\mathord{\mbox{\makebox[0pt][l]{\raisebox{.4ex}{$\uparrow$}}$\uparrow$}}}
\newcommand{\ua}{\mathord{\uparrow}}
\newcommand{\da}{\mathord{\downarrow}}
\newcommand{\rom}[1]{\rm{\uppercase\expandafter{\romannumeral #1}}}

\begin{document}

\begin{frontmatter}

%% Title,\ authors and addresses



\title{power structures of directed spaces\tnoteref{t1} }
\tnotetext[t1]{Research supported by NSF of China (Nos. 11871353, 12001385).}

\author{Xiaolin Xie}
\ead{xxldannyboy@163.com}

\author{Yuxu Chen\corref{cor}}
\ead{yuxuchen@stu.scu.edu.cn}
\cortext[cor]{Corresponding author}
\author{Hui Kou\corref{cor}}
\ead{kouhui@scu.edu.cn}
\address{Department of Mathematics,\ Sichuan University,\ Chengdu,610064, China }

%\address[add1]{Department of mathematics, Sichuan University, Chengdu,610064, China}
%\author[add1]{Yu Chen}
%\ead{eidolon_chenyu@126.com}

\begin{abstract}
Powerdomains in domain theory plays an important role in modeling the semantics of  nondeterministic  functional programming languages.\ In this paper,\ we extend the notion of  powerdomain to the category of  directed spaces,\ which is equivalent to the notion of the\ $T_0$\ monotone-determined space\ \cite{EN2009}.\ We define the notion of upper,\ lower and convex powerspace of a directed space by the way of free algebras.\ We show that the upper,\ lower and convex powerspace over any directed space exist and give their concrete structures.\ Generally,\ the upper,\ lower and convex powerspaces of a directed spaces are different from the upper,\ lower and convex powerdomains of a dcpos endowed with the Scott topology and the observationally-induced upper and lower powerspaces introduced by Battenfeld and Sch\"{o}der in 2015.\\\
\end{abstract}

\end{frontmatter}

\section{Introduction}

Powerdomain is one of the most important part of domain theory.\ Its purpose is to provide a mathematical model for the semantics of  nondeterministic  functional programming languages.\ There are three classical power structures in domain theory,\ they are the lower powerdomain(\ also known as the Hoaer powerdomain),\ the upper powerdomain(\ also known as the Smyth powerdomain)\ and the convex powerdomain(\ also known as the Plotikin powerdomain).\ Moreover,\ each power structure has a standard topological representation.\ In recent years,\ papers\ \cite{GK2017,YK2014-1,YK2014-2,YK2014-3}\ have made a lot of generalizations on these power structures.\ Generally speaking,\ these power structures are free algebras generated by domain respect to some binary operation.\ In 2015,\ I.\ Batenfeld and M.\ Sch\"{o}der\ (\cite{BS2015})introduced power structure in general topological spaces,\ and in this paper,\ Boolean algebra\ {\bf 2}\ is endowed with Sierpinski topology as an observable structure,\ and then the upper power structure and lower power structure are defined,\ which is called the observationally-induced upper(lower)\ powerspace.\ The work of paper\ \cite{BS2015}\ makes it possible for the  category of general topological space to be applied to express the nondeterministic semantics of functional programming languages.\ But the method of power structure induced by Boolean algebra\ {\bf 2}\ is complicated and difficult to understand.\ In this paper,\ we consider a class of special topological spaces,\ the directed space\ (equivalent to those in\ \cite{EN2009}).\ It is proved in\ \cite{YYK2015}\ that directed spaces contain the basic objects of domain theory,\ all directed complete posets endowed with the Scott topology,\ and continuous functions,\ which forms a cartesian closed category.\ Thus,\ the directed space is an extended framework of domain theory.\ In this paper,\ we are going to extend the structure of powerdomain to the category of directed spaces,\ and define the concept of directed upper (lower,\ convex) powerspace by means of free algebra.\ We show that the upper,\ lower and convex powerspace over any directed space exist and give their concrete structures.\ Generally,\ the upper,\ lower and convex powerspaces of a directed spaces are different from the upper,\ lower and convex powerdomains of a dcpos endowed with the Scott topology and the observationally-induced upper and lower powerspaces.\\

In 1991,\ Heckmann introduced an algebraic method and it does not rely on any explicit representations of the powerdomains\ \cite{HECK91}.\ In the last part of this article,\ we will discuss the commutativity of the directed upper and lower functors.\\

\section{Preliminaries}
Now,\ we introduce the concepts needed in this article.\ On\ domain\ theory,\ topology,\ and category theory,\ see\ \cite{AJ,  GHK,  Mac71}.\ Let\ $P$\ be a nonempty set.\ A relation\ $\leq$\  on\ $P$\ is called a partial order,\ if\ $\leq$\ satisfies reflexivity\ ($x\leq x$),\ transitivity\ ($x\leq y\ \&\ y\leq z\Rightarrow x\leq z$)\ and antisymmetry\ ($x\leq y\ \&\ y\leq x\Rightarrow x=y$).\ $P$\ is called a partial ordered set(poset) if $\ P\ $is endowed with some partial order\ $\leq$.\ Given\ $A\subseteq P$,\ denote\ $\da A=\{x\in P: \exists a\in A,\  \ x\leq a\}$,\ $\ua A=\{x\in P: \exists a\in A,\ a\leq x\}$.\ We say\ $A$\ is a lower set\ (upper set)\ if\ $A=\da A$ ($A=\ua A$).\ A nonempty set\ $D\subseteq P$\ is called a directed set if each finite nonempty subset of \ $D$\ has upper bound in\ $D$.\ Particularly,\ a poset is called a directed complete poset if each directed subset has a supremum(denoted by\ $\bigvee D$),\ abbreviated as\ dcpo.\ The subset\ $U$\ of poset\ $P$\ is called a\ Scott\ open set if\ $U$\ is an upper set and for each directed set\ $D\subseteq P$,\ which\ $\bigvee D$\ exists and belongs to\ $U$,\ then $U\cap D\not=\emptyset$.\ The set of all Scott open sets of poset\ $P$\ is a toplology on\ $P$,\ which is called the Scott topology and denoted by\  $\sigma(P)$.\ Suppose\ $P,\ E$\ are two posets,\ a fuction\ $f:P\longrightarrow E$\ is called\ Scott\ continuous if it is continuous respect to\ Scott\ topology\ $\sigma(P)$\ and\ $\sigma(E)$.

All topological spaces in this paper are required \ $T_0$\ separation.\ A net of a topological space\ $X$\ is a map\ $\xi: J\longrightarrow X$,\ here\ $J$\ is a directed set.\ Thus,\ each directed subset of a poset can be regarded as a net,\  and its index set is itself.\ Usually,\ we denote a net by\ $(x_j)_{j\in J}$\ or\ $(x_j)$.\ Let\ $x\in X$,\ saying\ $(x_j)$\ converges to\ $x$,\ denote by\ $(x_j)\rightarrow x$\ or\ $x\equiv \lim x_j$,\ if\ $(x_j)$\ is eventually in every open neighborhood of\ $x$,\ that is,\ for each given open neighborhood \ $U$\ of\ $x$,\ there exists\ $j_0\in J$\ such that for every\ $j\in J$,\ $j\geq j_0\Rightarrow x_j\in U$.

Let\ $X$\ be a\ $T_{0}$\ topological space,\ its topology is denoted by\ $O(X)$,\ the specialization order on\ $X$\ is defined as follows:
$$\forall x,\ y \in X,\ x\sqsubseteq y  \Leftrightarrow x\in \overline{\{y\}}$$
here,\ $\overline{\{y\}}$\ means the closure of\ $\{y\}$.
From now on,\ the order of a\ $T_{0}$\ topological space  always indicates the specialization order\ "$\sqsubseteq$".\ Here are some basic properties of specialization order.
\begin{proposition}{\rm\cite{AJ, GHK}}
	For a\ $T_{0}$\ topological space\ $X$,\ the following are always true:
	\begin{enumerate}
		\item For each open set\ $U\subseteq X,\ U=\ua U;$
		\item For each closed set\ $A\subseteq X,\ A=\da A;$
		\item Suppose\ $Y$\ is another\ $T_{0}$\ topological space,\ and\ $f:X\rightarrow Y$\ is a continuous function from\ $X$\ to\ $Y$.\ Then for each\ $x,\ y\in X$,\ $x\sqsubseteq y \Rightarrow f(x) \sqsubseteq f(y)$,\ that is every continuous fuction is monotone.
	\end{enumerate}
\end{proposition}

Suppose\ $X$\ is a\ $T_{0}$\ space,\ then every directed set\ $D\subseteq X$\ can be regarded as a net of\ $X$,\ we use\ $D\rightarrow x$\ or\ $x\equiv \lim D$\ to represent\ $D$\ converges to\ $x$.\ Define notation
\begin{center}
	$D(X)=\{(D,\   x):x\in X,\   D$\ is a directed subset of\ $X$\ and\ $D\rightarrow x \}$.
\end{center}
It is easy to verify that,\ for each\ $x,\   y\in X$,\ $x\sqsubseteq y\Leftrightarrow \{y\}\rightarrow x$.\ Therefore,\ if\ $x\sqsubseteq y$\ then\ $(\{y\},\   x)\in D(X)$.\ Next,\ we give the concept of directed space.

\begin{definition}{\rm\cite{YYK2015}}
	Let\ $X$\ be a\ $T_{0}$\ space.
	\begin{enumerate}
		\item[(1)] A subset\ $U$\ of\ $X$\ is called a directed open set if
		\ $\forall (D,\   x)\in D(X),\ x\in U \Rightarrow D\cap U\neq \emptyset$.\ Denote all directed open sets of \ $X$\ by\ $d(X)$.
		\item[(2)] \ $X$\ is called a directed space if each directed open set of\ $X$\ is an open set,\ that is,\ $d(X)=O(X)$.
	\end{enumerate}
\end{definition}

\begin{remark}\ \begin{enumerate}
		\item[(1)]\ Each open set of a\ $T_{0}$\ space is directed open,\ but the contrary is not necessarily true.\ For example,\ suppose\ $Y$\ is a non-discrete\ $T_1$\ topological space,\ its specialization order is diagonal,\ that is,\ $\forall x,\ y\in Y$,\ $x\sqsubseteq y\Leftrightarrow x=y$.\ Thus,\ all subsets of\ $Y$\ are directed open.\ We notice that\ $Y$\ is non-discrete,\ at least one directed open set is not an open set.
		\item[(2)] The definition of directed space here is equivalent to the $T_0$\ monotone determined space defined in\ {\rm \cite{EN2009}}.
		\item[(3)] Every poset endowed with the\ Scott\ topology is a directed space\ {\rm\cite{Kou15, YYK2015}},besides,each Alexandroff space is a directed space.\ Thus,\ the directed space extends the concept of the Scott topology.
	\end{enumerate}
\end{remark}

Next,\ we introduce the directed continuous function.

\begin{definition}\label{def dc}
	Suppose\ $X,\ Y$\ are two\ $T_{0}$\ spaces.\ A function\ $f:X\longrightarrow Y$\ is called directed continuous if it is monotone and preserves all limits of directed set of\ $X$;\ that is,\ $(D,\   x)\in D(X)\Rightarrow(f(D),\ f(x))\in D(Y)$.
\end{definition}

Here are  some characterizations of the directed continuous functions.
\begin{proposition}\label{prop dc}{\rm\cite{YYK2015}}
	Suppose\ $X,\ Y$\ are two\ $T_{0}$\ spaces.\ $f:X\longrightarrow Y$\ is a function between\ $X$\ and\ $Y$.
	\begin{enumerate}
		\item[(1)]$f$\ is directed continuous if and only if\ $\forall U\in d(Y),\ f^{-1}(U)\in d(X)$.
		\item[(2)]If\ $X,\ Y$\ are directed spaces,\ then \ $f$\ is continuous if and only if it is directed continuous.
	\end{enumerate}
\end{proposition}

Now we introduce the product of directed spaces.

Suppose\ $X,\ Y$\ are two directed spaces.\ Let\ $X\times Y$\ represents the cartesian product of\ $X$\ and\ $Y$,\ then we have a natural partial order on it:\ $\forall (x_{1},\ y_{1}),\ (x_{2},\ y_{2})\in X\times Y$,\ $$(x_{1},\ y_{1})\leq(x_{2},\ y_{2})\iff x_{1}\sqsubseteq  x_{2},\ y_{1}\sqsubseteq y_{2}.$$
which is called the pointwise order on\ $X\times Y$.\ Now,\ we define a topological space\ $X\otimes Y$\ as follows:
\begin{enumerate}
	\item The underlying set of\ $X\otimes Y$\ is\ $X\times Y$;
	\item The topology on\ $X\times Y$\ is generated as follows:\ for each given\ $\leq$-\ directed set\ $D\subseteq X\times Y$\ and\ $(x,\ y)\in X\times Y$,\ $$D\rightarrow(x,\ y)\in X\otimes Y\iff\pi_{1}D\rightarrow x\in X,\ \pi_{2}D\rightarrow y\in Y,\   $$
	That is,\ a subset\ $U\subseteq X\times Y$\ is open if and only if for every directed limit defined as above\ $D\rightarrow (x,\ y)$,\ $(x,\ y)\in U\Rightarrow U\cap D\neq\emptyset$.
\end{enumerate}
\begin{theorem}\label{theorem opc}{\rm\cite{YYK2015}}
	Suppose\ $X$\ and\ $Y$\ are two directed spaces.
	\begin{enumerate}
		\item The topological space\ $X\otimes Y$\ defined as above is a directed space and satisfies the following properties:\ the specialization order on\ $X\otimes Y$\ equals to the pointwise order on\ $X\times Y$,\ that is\ $\sqsubseteq =\leq$.
		
		\item Suppose\ $Z$\ is another directed space,\ then\ $f:X\otimes Y\longrightarrow Z$\ is continuous if and only if it is continuous in each variable separately.
	\end{enumerate}
\end{theorem}

Denote the category of all directed spaces with continuous functions by\ ${\bf Dtop}$.\ It is proved in\ \cite{YYK2014,YYK2015} that,\ ${\bf Dtop}$\ contains all posets endowed with the Scott topology and\ ${\bf Dtop}$\ is a\ cartesian closed category;\ specifically,\ the categorical products of two directed spaces\ $X$\ and\ $Y$\ are isomorphic to\ $X\otimes Y$.\ So,\ the directed space is an extended framework of\ Domain Theory.

Let\ $P$\ be a\ dcpo,\ and\ $x,\ y\in P$.\ We say\ $x$\ way below $y$,\ if for each given directed set\ $D\subseteq P$,\ $y\leq\bigvee D$\ implies there exists some\ $d\in D$\ such that\ $x\leq d$.\ We write\ $\dda x=\{a\in P:a\ll x\}$,\ $\dua x=\{a\in P:x\ll a\}$.

\begin{definition} A\ dcpo\ $P$\ is called a continuous domain if for each\ $x\in P$,\ $\dda x$\ is directed and\ $x=\bigvee \dda x$.
\end{definition}

\begin{theorem}{\rm\cite{GHK}}
	Suppose\ $P$\ is a continuous\ domain.\ The followings hold:
	\begin{enumerate}
		\item[(1)] $\forall x,\ y\in P$,\ $x\ll y\Rightarrow \exists z\in P,\ x\ll z\ll y$.
		\item[(2)] $\forall x\in P$,\ $\dua x$\ is a\ Scott\ open set.\ Particularly,\ $\{\dua x:x\in P\}$\ is a base of\ $(P,\ \sigma(P))$.
	\end{enumerate}
\end{theorem}
\vskip 3mm

\section{The directed lower powerspaces of directed spaces}
As mentioned above,\ directed space is a extended framework of\ domain\ theory,\ just like the work done in article\ \cite{BS2015},\ extending the powerdomain in the category of directed space is very meaningful.\ In this section,\ we will construct the directed lower powerspace of the directed space,\ which is a free algebra generated by the inflationary operation of the directed space.
	
\begin{definition}{\rm\cite{XXLK2020}}\label{def disl}
		Let\ $X$\ be a directed space.
\begin{enumerate}
\item[(1)] A binary operation\ $\oplus :X\otimes X\rightarrow X$\ on\ $X$\ is called an inflationary operation if it is continuous and satisfies the following four conditions: $\forall x,\ y,\ z\in X$,
		\begin{enumerate}
			\item[(a)] $x\oplus x=x$,\  			
\item[(b)] $(x\oplus y)\oplus z=x\oplus (y\oplus z)$,
\item[(c)] $x\oplus y=y\oplus x$,
\item[(d)] $x\oplus y\geq x$.
  		\end{enumerate}
\item[(2)] If\ $\oplus$\ is a inflationary operation on\ $X$,\ then\ $(X,\ \oplus)$\ is called a directed inflationary semilattice,\ that is,\ directed inflationary semilattices are those directed spaces with inflationary operations.
\end{enumerate}
	\end{definition}

By Theorem\ \ref{theorem opc}(2),\ the operation\ $\oplus$\ on a directed space\ $X$\ is continuous if and only if it is monotone and for each given\ $ x,\ y\in X$\ and directed set\ $D\subseteq X$,\ $\ x\equiv \lim D$\ implies\ $x\oplus y\equiv \lim (D\oplus y)$.\ Here,\ $D\oplus y=\{d\oplus y:d\in D\}$.

\vskip 3mm
Here are two examples.

\begin{example}{\rm\cite{XXLK2020}}
\begin{enumerate}
 \item[(1)] Suppose\ $P$\ is a poset endowed with\ the Scott\ topology,\ and for each\ $a,\ b\in P$,\ the supremum of\ $a$\ and\ $b$\ exists in\ $P$\ (denoted by\ $a\vee b$\ ).\ Then\ $(P,\ \vee)$\ is a directed inflationary semilattice.
 \item[(2)] Let\ $I=[0,\ 1]$\ (the unit interval),\ let\ ${\mathcal{A}}$\ denote the topology generated by\ $\{[a,\  1]:a\in I\}$,\ then\ $\mathcal{A}$\ is the Alexandroff topology on\ $I$.\ It is easy to check that,\ $(I,\ {\mathcal A})$\ is a directed space,\ and\ $(I,\ max)$\ is a directed inflationary semilattice  endowed with\ $\mathcal A$.
 \end{enumerate}
 \end{example}

\begin{definition}{\rm\cite{XXLK2020}}
Suppose\ $(X,\ \oplus),\ (Y,\ \uplus)$\ are two directed inflationary semilattices,\ $f:(X,\ \oplus)\rightarrow (Y,\ \uplus)$\ is called an inflationary homomorphism between\ $X$\ and\ $Y$,\ if\ $f$\ is continuous and\ $f(x\oplus y)=f(x)\uplus f(y)$\ holds,\ $\forall x,\ y\in X$.
\end{definition}

Denote the category of all directed inflationary semilattices and inflationary homomorphisms by\ ${\bf Disl}$.\ Then\ ${\bf Disl}$\ is a subcategory of \ ${\bf Dtop}$.

\begin{lemma}{\rm\cite{XXLK2020}}
Suppose\ $(X,\ \oplus)$\ is a directed inflationary semilattice,\ then\ $\oplus =\vee_{\sqsubseteq} $.\ Here,\ $\forall x,\ y \in X,\ x\vee_{\sqsubseteq}y$\ means the supremum of\ $x$\ and\ $y$\ respect to the specialization order\ $\sqsubseteq$\ on\ $X$\ (calling it sup operation).\ Conversely,\ suppose\ $X$\ is a directed space such that for each\ $x,\ y \in X $,\ $x\vee_{\sqsubseteq}y$\ exists,\ and the continuity of \ $\vee_{\sqsubseteq} $\ will naturally imply that\ $(X,\ \vee_{\sqsubseteq})$\ is a directed inflationary semilattice.
\end{lemma}
\noindent{\bf Proof}  By definition\ \ref{def disl} ,$\forall x,\ y\in X,\ x\oplus y\geq x,\ y$,\ that is\ $x\oplus y$\ is a upper bound of\ $\{x$,\ $y\}$.\ Suppose\ $z$\ is another arbitrary upper bound of\ $\{x,\ y\}$,\ By Theorem\ \ref{theorem opc}(1),\ the pointwise order equals to the specialization order of\ $X\otimes X$,then\ $(x,\ y)\sqsubseteq (z,\ z)$.\ By the continuity and idempotence of the inflationary operation,\ we have\ $x\oplus y\sqsubseteq z\oplus z=z$.\ That is\ $x\oplus y$\ is the supremum of\ $\{x,\ y\}$,\ which means\ $x\oplus y=x\vee_{\sqsubseteq}y$.\ Conversely,\ continuous sup operation naturally satisfy all conditions in definition\ \ref{def ddsl},\ thus,\ we get the conclusion.\ $\Box$

\vskip 3mm
The above results show that a directed inflationary semilattice\ $(X,\ \oplus)$\ is just a directed space with a continuous sup operation satisfy\ $\oplus=\vee_{\sqsubseteq}$.\ Since the order on the directed space in this paper is always the specialization order,\ we will use the symbol\ $\vee$\ instead of\ $\vee_{\sqsubseteq}$\ in the following part.\ Therefore,\ a directed inflationary semilattice is always represented by a tuple of the form\ $(X,\ \vee)$,\ here\ $X$\ is a directed space,\ $\vee$\ represents the continuous sup operation on\ $X$.
\vskip 3mm
Next,\ we give the definition of directed lower powerspace.

\begin{definition}{\rm\cite{XXLK2020}}\label{def dlps}
  Suppose\ $X$\ is a directed space.\ A directed space\ $Z$\ is called the directed lower powerspace of\ $X$\ if and only if the following two conditions are satisfied:
\begin{enumerate}
\item[(1)] $Z$\ is a directed inflationary semilattice,\ that is the sup operation\ $\vee$\ on\ $Z$\ exists and which is continuous,
\item[(2)] There is a continuous function\ $i:X\longrightarrow Z$\ satisfying:\ for an arbitrary directed inflationary semilattice\ $(Y,\ \vee)$\ and continuous function\ $f:X\longrightarrow Y$,\ there exists an unique inflationary homomorphism\ $\bar{f}:(Z,\ \vee)\rightarrow (Y,\ \vee)$\ such that\ $f=\bar{f}\circ i$.
\end{enumerate}
\end{definition}

By the definition above,\ if directed inflationary semilattices\ $(Z_1,\ \vee)$\ and\ $(Z_2,\ \vee)$\ are both the directed lower powerspaces of\ $X$,\ then there exists a topological homomorphism which is also a inflationary homomorphism\ $g:Z_1\rightarrow Z_2$.\ Therefore,\ up to of order isomorphism and topological homomorphism,\ the directed lower powerspace of a directed space is unique.\ Particularly,\ we denote the directed lower powerspace of each directed space\ $X$\ by\ $P_L(X)$.

\vskip 3mm

Now,\ we will prove the existence of the directed lower powerspace of each directed space\ $X$\ by way of concrete construction.

\vskip 3mm

	Let\ $X$\ be a directed space.\ Set
$$LX=\{\da F :F\subseteq _{fin} X\},\  $$
here,\ $F\subseteq _{fin} X$\ is an arbitrary nonempty finite subset of \ $X$.\ Define an order\ $\leq_L$\ on\ $LX$\ as follows:
$$\da F_{1}\leq_L \da F_{2}
	\iff \da F_{1}\subseteq \da F_{2}.$$
	Let\ $\mathcal{D}\subseteq LX$\ be a directed set\ (respect to order\ $\leq_L$),\ $\da F\in LX$.\ Define a convergence notation\ $\mathcal{D}\Rightarrow_L \da F$\ as follows:
	\begin{center}
		$\mathcal{D}\Rightarrow _{L}\da F \iff \forall a \in  F,\ $\ there exists a directed set\ $D_{a}$\ of\ $X$\ satisfying $D_{a}\subseteq \bigcup \mathcal{D}$\ and\ $D_{a}\rightarrow a$.
	\end{center}
A subset\ $\mathcal{U}\subseteq LX$\ is called a\ $\Rightarrow_L$\ convergence open set of\ $LX$\ if and only if for each directed subset\ $\mathcal{D}$\ of\ $LX$\ and\ $\da F\in LX$,\ \ $\mathcal{D}\Rightarrow_{L}\da F\in \mathcal{U}$\ implies\ $\mathcal{D}\cap \mathcal{U}\neq \emptyset$.\ Denote all\ $\Rightarrow_L$\ convergence open sets of\ $LX$\ by\ $O_{\Rightarrow_{L}}(LX)$.

\begin{proposition}{\rm\cite{XXLK2020}}\label{prop LX}
	Suppose\ $X$\ is a directed space,\ the following are true:
	\begin{enumerate}
		\item[(1)] $(LX,\ O_{\Rightarrow_L}(LX))$\ is a topological space,\ abbreviated as\ $LX$.
		\item[(2)] The specialization order\ $\sqsubseteq $\ on\ $(LX,\ O_{\Rightarrow_L}(LX))$\ equals to\ $\leq_L$.
		\item[(3)] $(LX,\ O_{\Rightarrow_L}(LX))$\ is a directed space,\ that is\ $O_{\Rightarrow_L}(LX)=d(LX)$.
		
	\end{enumerate}
\end{proposition}
\noindent{\bf Proof}
(1)  Firstly\ $\emptyset,\ LX\in O_{\Rightarrow_L}(LX)$.\ Suppose\ $\mathcal{U}\in O_{\Rightarrow_L}(LX)$,\ $\da F_1,\ \da F_2\in LX$.\ If\ $\da F_{1}\leq \da F_{2}$,\ then\ $\{\da F_{2}\}\Rightarrow_{L}\da F_{1}$\ is obviously hold.\ Therefore\ $\da F_{1}\in\mathcal{U}$\ will imply\ $\da
F_{2}\in\mathcal{U}$.\ This means\ $\mathcal{U}$\ is an upper set respect to order\ $\leq_L$\ on\ $LX$.\ Suppose we have\ $\mathcal{U}_{1},\  \mathcal{U}_{2}\in O_{\Rightarrow_L}(LX)$,\ and\ $\mathcal{D}$\ is a directed set in\ $LX$,\ $\da F\in LX$\ and\ $\mathcal{D} \Rightarrow \da F\in \mathcal{U}_{1}\cap\mathcal{U}_{2}$,\ there will exist\ $\da G_1\in \mathcal{D}\cap \mathcal{U}_{1}$\ and\ $\da G_2\in \mathcal{D}\cap \mathcal{U}_{2}$.\ Since\ $\mathcal{D}$\ is directed,\ there exists\ $\da G\in \mathcal{D}$\ such that\ $\da G_1,\  \da G_2\leq_L \da G$,\ taht is\ $\da G \in \mathcal{D}\cap \mathcal{U}_{1}\cap\mathcal{U}_{2}$.\ By the same way,\ $O_{\Rightarrow_L}(LX)$\ is closed under the arbitrary union.\ Thus,\ $O_{\Rightarrow_L}(LX)$\ is a topology on\ $LX$.

\vskip 3mm

(2) Let\ $\da F_1,\ \da F_2\in LX$.\ If\ $\da F_{1}\leq_L \da F_{2}$,\ according the proof of\ (1),\ every\ $\Rightarrow_L$\ convergence open set is an upper set respect to partial order\ $\leq_L$,\ then\ $\da F_1\in \overline{\{\da F_2\}}$,\ that is,\ $\da
F_{1}\sqsubseteq \da F_{2}$.\
On the other hand,\ let\ $\da F_{1}\sqsubseteq \da F_{2}$.\ We want to prove\ $\da F_{1}\leq_L \da F_{2}$,\ that is\ $\{\da F\in LX: \da F\subseteq \da F_2\}$\ is a closed set of\ $LX$\ respect to the topology\ $O_{\Rightarrow_L}(LX)$,\ since\ $\{\da F\in LX: \da F\subseteq \da F_2\}\subseteq \overline{\{\da F_2\}}$,\ $\{\da F\in LX: \da F\subseteq \da F_2\}$\ is a closed set of\ $LX$\ respect to the topology\ $O_{\Rightarrow_L}(LX)$\ will imply\ $\{\da F\in LX: \da F\subseteq \da F_2\}=\overline{\{\da F_2\}}$,\ that is\ $\da F_{1}\leq_L \da F_{2}$.\ We are going to prove that\ $\{\da F\in LX: \da F\subseteq \da F_2\}$\ is closed respcet to topology\ $O_{\Rightarrow_L}(LX)$,\ that is to show\ $\mathcal{U}=LX\setminus\{\da F\in LX: \da F\subseteq \da F_2\}$\ is a\ $\Rightarrow_L$\ convergence open set in\ $LX$.\ Let\ ${\mathcal D}\subseteq LX$\ be a\ $\leq_L$\ -\ directed set and\ $\mathcal{D}\Rightarrow_L \da G\in\mathcal{U}$.\ Let\ $G=\{a_1,\  a_2,\  \ldots,\  a_k\}$.\ By the definition of\ $\Rightarrow_L$\ convergence,\ we have finitary directed subsets\ $D_i\subseteq\bigcup {\mathcal D}$\ in\ $X$\ such that\ $D_i\rightarrow a_i$,\ $i=1,\  2,\  \ldots,\  k$.\ By contradiction,\ we suppose that\ ${\mathcal D}\cap {\mathcal U}=\emptyset$.\ Then\ $\bigcup {\mathcal D}\subseteq \da F_2$,\ that is\ $D_i\subseteq \da F_2$,\  $i=1,\  2,\  \ldots,\  k$.\ Since\ $\da F_2$\ is closed in\ $X$,\ the limit points of\ $D_i$\ are in\ $\da F_2$,\ that is\ $G=\{a_1,\  a_2,\  \ldots,\  a_k\}\subseteq \da F_2$,\ and,\ $\da G\leq_L\da F_2$.\ This contradicts with\ $\da G\in {\mathcal U}$.\ Therefore,\ $\mathcal{U}=LX\setminus\{\da F\in LX: \da F\subseteq \da F_2\}$\ is a\ $\Rightarrow_L$\ convergence open set in\ $LX$.

\vskip 3mm

(3) For an arbitrary topological space\ $X$,\ we have\ $O(X)\subseteq d(X)$,\ then\ $O_{\Rightarrow_L}(LX)\subseteq d(LX)$.\ On the other hand,\ according the definition of\ $\Rightarrow_L$\ convergence,\ if directed set\ $\mathcal{D}\Rightarrow_{L}\da F$\ in\ $LX$,\ according to\ (2),\ $\mathcal{D}$\ convergents to\ $\da F$\ respect to the topology\ $O_{\Rightarrow_L}(LX)$.\ Then,\ by the definition of directed open set,\ $\mathcal{D}\Rightarrow_{L}\da F \in\mathcal{U}\in d(LX)$,\ will imply\ $\mathcal{U}\cap \mathcal{D}\neq \emptyset$.\ This means\ $\mathcal{U}\in O_{\Rightarrow_L}(LX)$,\ thus\ $O_{\Rightarrow_L}(LX)= d(LX)$,\ that is,\ $(LX,\  O_{\Rightarrow_L}(LX))$ is a directed space.\ $\Box$

\vskip 3mm

\begin{proposition}{\rm\cite{XXLK2020}}\label{prop LXdc}
Suppose\ $X,\ Y$\ are two directed spaces.\ Then function\ $f:LX\rightarrow Y$\ is continuous if and only if for each directed set\ ${\mathcal D}\subseteq LX$\ and\ $\da F\in LX$,\ ${\mathcal D}\Rightarrow_L\da F$\ implies\ $f(\mathcal{D})\rightarrow f(\da F)$.
\end{proposition}
\noindent{\bf Proof} Since\ $\Rightarrow_L$\ convergence will lead to \ $O_{\Rightarrow_L}(LX)$\ topological convergence,\ the necessity is obviously.\ we are going to prove the sufficiency.\ First to check that\ $f$\ is monotone.\ If\ $\da F_1,\  \da F_2\in LX$\ and\ $\da F_1\leq_L\da F_2$,\ then\ $\{\da F_2\}\Rightarrow_L\da F_1$,\ by the given condition,\ $\{f(\da F_2)\}\rightarrow f(\da F_1)$,\ thus\ $f(\da F_2)\sqsubseteq f(\da F_1)$.\ Suppose\ $U$\ is an open set of\ $Y$\ and the directed set\ ${\mathcal D}\Rightarrow_L\da F\in f^{-1}(U)$,\ then\ $f(\mathcal{D})$\ is a directed set of\ $Y$\ and\ $f(\mathcal{D})\rightarrow f(\da F)\in U$,\ there exists a$ \da F\in {\mathcal D}$\ such that\ $f(\da F)\in U$.\ Thus,\ $\da F\in {\mathcal D}\cap f^{-1}(U)$.\ According to the definition of\ $\Rightarrow_L$\ convergence open set,\ $f^{-1}(U)\in O_{\Rightarrow_L}(LX)$,\ that is\ $f$\ is continuous.\ $\Box$

\begin{theorem}{\rm\cite{XXLK2020}}
Let\ $X$\ be a directed space.\ Then\ $(LX,\ O_{\Rightarrow_L}(LX))$\ respect to the set union operation\ $\cup $\ is a directed inflationary semilattice.
\end{theorem}
\noindent{\bf Proof} According to Proposition\ \ref{prop LX},\ $(LX,\  O_{\Rightarrow_L}(LX))$\ is a directed space.\ we will prove that\ $\cup $\ is an inflationary operation.\ For arbitrary\ $\da F_1,\ \da
F_2\in LX$,\ then\ $\da F_1\cup \da F_2=\da (F_1\cup F_2)\in LX$.\ Obviously,\ $\cup$\ satisfy the  conditions \ (a),\  (b),\  (c),\  (d)\ in Definition\ \ref{def disl},\ we now prove the continuity of \ $\cup$.\ The monotonicity of\ $\cup$\ is evident.\ By Theorem\ \ref{theorem opc}(2)\ and Proposition\ \ref{prop LXdc},\ we only need to prove that,\ For each directed set\ ${\mathcal D}\subseteq LX$\ and\ $\da F,\  \da G\in LX$,\ \ ${\mathcal D}\Rightarrow_L\da F$,\ will imply\ $G\cup {\mathcal D}\Rightarrow_L(\da G\cup \da F)=\da(G\cup F)$.\ Here,\ $G\cup {\mathcal D}=\{\da (G\cup F'):\da F'\in \mathcal D\}$\ is a directed set.\ For each\ $a\in G\cup F$,\ if\ $a\in G$,\ then\ $\{a\}\rightarrow a$;\ if\ $a\in F$,\ since\ ${\mathcal D}\Rightarrow_L\da F$,there exists\ $D\subseteq \bigcup{\mathcal D}$\ such that\ $D\rightarrow a$.\ According to the definition of\ $\Rightarrow_L$\ convergence,\ we have\ $G\cup {\mathcal D}\Rightarrow_L\da(G\cup F)$.\ $\Box$

\begin{remark}
We can directly check that for arbitrary\ $\da F_1,\ \da F_2\in LX,\ \da F_1\cup\da F_2=\da (F_1\cup F_2)$\ is well-defined although each\ $\da F$\ may generated by different\ $F$.
\end{remark}

The following theorem is the main result in this section.

\begin{theorem}\label{theorem LX}\cite{XXLK2020}
	Suppose\ $X$\ is a directed space,\ then\ $(LX,\ O_{\Rightarrow_L}(LX))$\ is the lower powerspace of\ $X$,\ that is,\ endowed with topology\ $O_{\Rightarrow_L}(LX)$,\ $(LX,\ \cup)\cong P_L(X)$.	
\end{theorem}
\noindent{\bf Proof} Define function\ $i:X\rightarrow LX$,\ $\forall x\in X$,\ $i(x)=\da x$.\ We prove the continuity of\ $i$.\ It is evident that\ $i$\ is monotone.\ If\ $D\subseteq X$\ and\ $x\in X$\ satisfy\ $D\rightarrow x$.\ Let\ $\mathcal{D}=\{\da d:d\in D\}$,\ then\ $\mathcal{D}$\ is a directed set in\ $LX$\ and\ ${\mathcal D}\Rightarrow_{L}\da x$.\ But\ $i(D)=\mathcal{D}$,\ so\ $i(D)\Rightarrow_{L}\da x=i(x)$.\ According to Proposition\ \ref{prop LXdc},\ $i$\ is continuous.

Let\ $(Y,\ \vee)$\ be an arbitrary directed inflationary semilattice,\ $f:X\rightarrow Y$\ is a continuous function.\ Define\ $\bar{f}:LX\rightarrow Y$\ as follows:\ $\forall\  \da F\in LX$\ (let\ $F=\{a_1,\  a_2,\  \ldots,\  a_n\}$\ ),\ $$\bar{f}(\da F)=f(a_1)\vee f(a_2)\vee\cdots\vee f(a_n)=\bigvee\limits_{a\in F}f(a).$$
$\bar{f}$\ obvious well-defined.\ Especially,\ denote\ $\bar{f}(\da F)=\vee f(F)$.

(1) $f=\bar{f}\circ i$.

For arbitrary\ $x\in X$,\ $(\bar{f}\circ i)(x)=\bar{f}(i(x))=\bar{f}(\da x)=f(x)$.

(2) $\bar{f}$\ is an inflationary homomorphism,\ that is,\ $\bar{f}$\ is continuous and for arbitrary\ $\da F_1,\  \da F_2\in LX$,\  $\bar{f}(\da F_1\cup \da F_2)=\bar{f}(\da F_1)\vee \bar{f}(\da F_2)$.

First,we prove that\ $\bar{f}$\ preserves the union operation.\ Let\ $\da F_1,\ \da F_2\in LX$.\ Then\ $\bar{f}(\da F_1\cup\da F_2)=\bar{f}(\da(F_1\cup F_2))=\vee f(F_1\cup F_2)=(\vee f(F_1)\vee(\vee f(F_2))=\bar{f}(\da F_1)\vee\bar{f}(\da F_2)$.\ Next,\ we prove the continuity of\ $\bar{f}$.\ Since\ $\vee$\ is the sup operation,\ $\bar{f}$\ is monotone.\ Suppose\ ${\mathcal D}\subseteq LX$\ and\ $\da G\in LX$\ satisfy\ ${\mathcal D}\Rightarrow_L\da G$.\ Let\ $G=\{b_1,\  b_2,\  \ldots,\  b_k\}$.\ By the definition of\ $\Rightarrow_L$\ convergence,\ for each\ $b_i\in G$,\ there exists a directed set\ $D_i\subseteq \bigcup{\mathcal D}$\ such that\ $D_i\rightarrow b_i$.\ By the continuity of\ $f$,\ we have\ $f(D_i)\rightarrow f(b_i)$,\  $i=1,\  2,\  \ldots,\  k$.\ since\ $\vee$\ is continuous,the following convergence hold in Y:\
$$f(D_{1})\vee f(D_2)\vee\cdots \vee f(D_{k})\rightarrow f(b_{1})\vee \dots \vee f(b_k).\ \\ (\ast)$$
Here,\ $$f(D_{1})\vee f(D_2)\cdots \vee f(D_{k})=\{f(d_{1})\vee f(d_2)\vee \cdots \vee f(d_k):(d_1,\  d_2,\  \ldots,\  d_k)\in \prod\limits_{i=1}^{k}D_i\}.\ $$
Suppose\ $(d_{1},\  d_2,\  \dots,\  d_k)
\in \prod\limits_{i=1}^{k}D_{i}$.\ For arbitrary\ $i\in\{1,\  2,\  \ldots,\  k\}$,\ there exists\ $\da F_{i}\in \mathcal{D}$\ such that\ $d_{i}\in \da F_{i}$.\ Since\ $\mathcal{D}$\ is directed,\ there exists some\ $\da F\in \mathcal{D}$ such that\ $\da F_{i}\subseteq  \da F$.\ Therefore\ $f(d_{1})\vee \dots \vee f(d_{n})\sqsubseteq
(\vee f(F_{1}))\vee \dots \vee (\vee f(F_{k}))\sqsubseteq \vee f(F)$.\ This means that,\ $$f(D_{1})\vee \dots \vee f(D_{n})\subseteq \da \{\vee f(F):\da F\in \mathcal{D}\}.\ \\ (\ast\ast)$$
Let\ $U\subseteq Y$\ be an open neighborhood of\ $f(b_{1}) \vee f(b_2)\vee\cdots \vee f(b_{k})$.\ By\ $(\ast)$,\ there exists some\ $(d_{1},\  d_2\ldots,\  d_{k})
\in \prod\limits_{i=1}^{k}D_{i}$\ such that\ $f(d_{1})\vee f(d_2)\vee
\cdots \vee f(d_{k})\in U$.\ For each open set is an upper set,\ by\ $(\ast\ast)$,\ there exists some\ $\da F\in {\mathcal D}$\ such that\ $\vee f(F)\in U$.\ That is\ $\bar{f}(\mathcal{D})\rightarrow \bar{f}(\da G) $.\ Therefore,\ by Proposition\ \ref{prop LXdc},\ $\bar{f}$\ is continuous.

(3) Inflationary homomorphism\ $\bar{f}$\ is unique.

Suppose there is an inflationary homomorphism\ $g:(LX,\  \cup)\rightarrow (Y,\ \vee)$\ satisfy\ $f=g\circ i$.\ Then\ $g(\da x)=f(x)=\bar{f}(\da x)$.\ For each\ $ \da F\in LX$\ (let\ $ F=(a_{1},\  \dots ,\  a_{n})$\ ),\ \begin{eqnarray*}
	g(\da F)&=& g(\da a_1\cup\da a_2\cdots\cup\da a_n)\\ &=&g(\da a_{1})\vee g(\da a_2)\vee \cdots \vee g(\da a_{n})\\ &=&\bar{f}(\da a_{1})\vee \bar{f}(\da a_2)\vee\cdots\vee\bar{f}(\da a_{n})\\ &=&\bar{f}(\da a_1\cup\da a_2\cdots\cup\da a_n)\\ &=&\bar{f}(\da F).
\end{eqnarray*}
Thus\ $\bar{f}$\ is unique.

In conclusion,\ according to Definition\ \ref{def dlps},\ endowed with topology\ $O_{\Rightarrow_L}(LX)$,\ the directed inflationary semilattice\ $(LX,\ \cup)$\ is the lower powerspace of\ $X$,\ that is,\ $P_L(X)\cong (LX,\ \cup)$.\ $\Box$
\vskip 3mm

  The lower powerspace is unique in the sense of order isomorphism and topological homomorphism,\ so we can directly denote the lower powerspace of each directed space\ $X$\  by\ $P_L(X)=(LX,\ \cup)$.

 Suppose\ $X,\ Y$\ are two directed spaces,\ $f:X\rightarrow Y$\ is a continuous function.\ Define map\ $P_L(f):P_L(X)\rightarrow P_L(Y)$\ as follows:\ $\forall \da F\in LX$,$$P_L(f)(\da F)=\da f(F).$$
 it is evident that\ $P_L(f)$\ is well-defined and order preserving.\ According to\ \ref{theorem LX},\ it is easy to check that\ $P_L(f)$\ is an inflationary homomorphism between these two lower powerspaces.\ If\ $id_X$\ is the identity function and\ $g:Y\rightarrow Z$\ is an arbitrary continuous function from\ $Y$\ to a directed space\ $Z$,\ then,\ $P_L(id_X)=id_{P_L(X)},\ P_L(g\circ f)=P_L(g)\circ P_L(f)$.\ Thus,\ $P_L:{\bf Dtop}\rightarrow {\bf Disl}$\ is a functor from\ ${\bf Dtop}$\ to\ ${\bf Disl}$.\ Let\ $U:{\bf Disl}\rightarrow {\bf Dtop}$\ be the forgetful functor,\ by Theorem\ \ref{theorem LX},\ we have the following result.

\vskip 3mm

 \begin{corollary}{\rm\cite{XXLK2020}}
$P_L$\ is a left adjoint of the forgetful functor\ $U$,\ that is,\ ${\bf Disl}$\ is a reflective subcategory of \ ${\bf Dtop}$.
 \end{corollary}
\vskip 3mm

\section{Relations Between Lower Powerspaces}
In this section,\ we will discuss the relation between the lower powerdomain of dcpo,\ observationally-induced lower powerspaces and directed lower powerspaces.

Suppose\ $X$\ is a topological space (denote the topology by\ $O(X)$\ ).\ Let\ $C(X)$\ be the set of all nonempty closed sets of\ $X$.\ Obviously,\ $C(X)$\ is closed under the set operation\ $\cup$.\ For each\ $U\in O(X)$,\ let
$$\langle U\rangle=\{A\in C(X):A\cap U\not=\emptyset\}.$$
It is easy to check that,\ $\{\langle U\rangle: U\in O(X)\}$\ is closed under finitary intersection and arbitrary union,\ it forms a topology on\ $C(X)$,\ called the lower Vietoris topology of\ $C(X)$\ and denoted by\ $V_L(C(X))$.\ Besides,\ $C(X)$\ is a\ dcpo\ respect to the inclusion order of set.\ Denote the\ Scott\ topology on\ $C(X)$\ by\ $\sigma(C(X))$.\ We can check that,\ the lower\ Vietoris\ topology\ $V_L(C(X))$\ is contained in the\ Scott\ topology\ $\sigma(C(X))$,\ that is,\ $V_L(C(X))\subseteq \sigma(C(X))$.\ Particularly,\ write\ $P_{OL}(X)=(C(X),\ \cup)$\ endowed with the lower\ Vietoris\ topology\ $V_L(C(X))$,\ write\ $H(X)=(C(X),\ \cup)$\ endowed with the Scott topology\ $\sigma(C(X))$.

\begin{theorem}{\rm \cite[Corollary IV-8.6]{GHK}}
  Let\ $P$\ is a\ dcpo\ endowed with the\ Scott\ topology\ $\sigma(P)$.\ Then\ $H(P)$\ isomorphic to the lower powerdomian of\ $P$\  (which is also called\ Hoare\ powerdomain\ ),\ that is,\ $H(P)$\ is free \ dcpo\ sup semilattice of generated by\ $P$.
\end{theorem}

\begin{theorem}{\rm \cite[Theorem 3.8]{BS2015}}
Let\ $X$\ be a topological space.\ Then\ $P_{OL}(X)$\ is isomorphic to the observationally-induced lower powerspace over\ $X$.
\end{theorem}

Let\ $P$\ be a dcpo,\ endowed with the\ Scott\ topology,\ then\ $(P,\ \sigma(P))$\ is also a directed space.\ By Theorem\ \ref{theorem LX},\ the directed lower powerspace\ $P_L(P)=(LP,\ \cup)$,\ endowed with topology\ $O_{\Rightarrow_L}(LP)$,\ compared with the lower powerdomain and the observationally-induced lower powerspace over\ $P$,\ in which the operation are both the union of sets\ $\cup$,\ $LP$\ is a subset of\ $C(P)$,\ and in general,\ $LP\not=C(P)$.\ For the topological structure,\ $P_L(P)$\ is neither necessarily a subspace of lower powerdomain\ $H(P)$,\ nor a subspace of the observationally-induced lower powerspace\ $P_{OL}(P)$.

The following example tells us that,\ in general,\ a\ dcpo\ $P$,\ endowed with\ Scott\ topology,\ its lower powerdomain is not equal to its observationally-induced lower powerspace,\ that is,\ $H(P)\not=P_{OL}(P)$.

\begin{example}{\rm\cite{XXLK2020}} Let\ $P=(\mathbb{N}\times\mathbb{N})\cup\{\infty\}$,\ here\ $\mathbb{N}$\ is the set of natural numbers.\ Define a partial order on\ $P$\ as follows:\ $\forall x,\ y\in P$,\  $x\leq y$\ if and only if one of the following conditions is true:
\begin{itemize}
\item $y=\infty$,
\item $\exists n_0\in\mathbb{N},\ x=(m,\ n_0),\ y=(m',\ n_0)$\ and\ $m'-m\geq 0$.

\end{itemize}
Evidently,\ $P$\ is a\ dcpo.\ And for each nonempty subset\ $A$\ ,\ $A$\ is a\ Scott\ closed set of \ $P$\ if and only if\ $A=P$\ or there exists an antichain\ $B$\ of\ $P$\ such that\ $A=\da B$.\ Then,\ we have\ $(\mathbb {N}\times \{n\})\cup\{\infty\}$\ is in\ $ C(P)$\ but not in \ $LP$,\ that is\ $LP\not=C(P)$.\ For each\ $n\in\mathbb{N}$,\ let\ $D_n=\mathbb{N}\times\{n\}$,\ then\ $D_n$\ is a directed set and\ $\bigvee(\mathbb{N}\times\{n\})=\infty$.\ For each\ $U\subseteq P$,\ $U$\ is a\ Scott\ open subset if and only if\ $U=\emptyset$\ or\ $U\cap D_n\not=\emptyset$.\ So,\ each nonempty\ Scott\ open set of\ $P$\ equals to the upper  set of its minimal elements\ ${\rm min}U$.\ Let\ $\mathcal{U}\subseteq C(P)$,\ then\ $\mathcal{U}$\ is a open set respect to the lower\ Vietoris\ topology if and only if there exists some\ $U\in \sigma(P)$\ such that$\ \mathcal{U}=\{B\in C(P):B\cap {\rm min}U\not=\emptyset\}$.\ Let
$$\mathcal{V}=\{A\in C(P):\exists n\in\mathbb{N},\ (5,  n)\in A\}\cup \{A\in C(P): \exists n\in\mathbb{N},\ (4,  n),  (4,  n+1)\in A\}.$$
we can easily check that,\ $\mathcal{V}$\ is a\ Scott\ open set of\ $C(P)$,\ and for arbitrary\ $V\in \sigma(P)$,\ $\mathcal{V}\not=\langle V\rangle$.\ Thus,\ $\mathcal{V}$\ is not an open set in lower\ Vietoris\ topology.\ That is,\ $H(P)\not=P_{OL}(P)$.
\end{example}
\vskip 3mm
Next,\ we are going to consider the relation between the lower powerdomains and the directed lower powerspaces.

Let\ $P$\ be a\ dcpo\ endowed with Scott\ topology\ $\sigma(P)$.\ Write
$$\sigma(C(P))|_{LP}=\{\mathcal{U}\cap LP: \mathcal{U}\in \sigma(C(X))\}$$
to denote the relative topology on\ $LP$\ from the\ Scott\ topology on\ $C(P)$.\ For arbitrary\ $\mathcal{V}\in O_{\Rightarrow_L}(LP)$,\ let\ $$\uparrow_{C(P)}\mathcal{V}=\{A\in C(P): \exists \da F\in \mathcal{V},\ \da F\subseteq A\}.$$
\vskip 3mm
\begin{proposition}{\rm\cite{XXLK2020}}\label{prop rstLX}
Let\ $P$\ be a\ dcpo\ endowed with\ Scott\ topology\ $\sigma(P)$.\ Then\ $\sigma(C(P))|_{LP}\subseteq O_{\Rightarrow_L}(LP)$.\ Particularly,\ $\sigma(C(P))|_{LP}=O_{\Rightarrow_L}(LP)$\ if and only if for each\ $\mathcal{V}\in O_{\Rightarrow_L}(LP)$,\ $\uparrow_{C(P)}\mathcal{V}\in\sigma(C(P))$.
\end{proposition}
\noindent{\bf Proof}
Let\ $\mathcal{U}\in \sigma(C(P))$,\ and a directed set\ ${\mathcal D}\subseteq LP$\ satisfy\ ${\mathcal D}\Rightarrow_L\da F\in \mathcal{U}\cap LP$.\ Then for each\ $a\in F$,\ there exists a directed set\ $D_a\subseteq \bigcup\mathcal{D}$\ such that\ $D_a\rightarrow a$.\ Thus,\ $\da F\subseteq \overline{\bigcup{\mathcal D}}$,\ Here\ $\overline{\bigcup{\mathcal D}}$\ means the\ Scott\ closure of\ $\bigcup{\mathcal D}$.\ But\ $\overline{\bigcup{\mathcal D}}$\ is just the supremum of\ $\bigcup{\mathcal D}$\ in\ $C(P)$,\ Then\ $\mathcal{D}\cap\mathcal{U}\cap LP=\mathcal{D}\cap \mathcal{U}\not=\emptyset$.\ It follows that\ $\mathcal{U}\cap LP\in O_{\Rightarrow_L}(LP)$,\ that is,\ $\sigma(C(P))|_{LP}\subseteq O_{\Rightarrow_L}(LP)$.

On the other hand,\ for arbitrary\ Scott\ closed set\ $A$\ with\ $A\in \mathcal{U}$,\ here\ $\mathcal{U}$\ is a\ Scott\ open set of\ $C(P)$.\ Let\ ${\mathcal F}(A)=\{\da F: F\subseteq_{fin}A\}$,\ then\ ${\mathcal F}(A)$\ is a directed set of\ $LP$\ and\ $A=\bigcup {\mathcal F}(A)$.\ We have a nonempty finite set\ $F$\ of\ $A$\ such that\ $\da F\in \mathcal{U}$.\ It means that $\mathcal{U}=\ua_{C(P)}(\mathcal{U}\cap LP)$.\ Thus,\ $\sigma(C(P))|_{LP}=O_{\Rightarrow_L}(LP)$\ if and only if\ $\forall \mathcal{V}\in O_{\Rightarrow_L}(LP)$,\ $\uparrow_{C(P)}\mathcal{V}\in\sigma(C(P))$.\ $\Box$

\vskip 3mm

Suppose\ $A$\ is an arbitrary subset of\ dcpo\ $P$.\ Let\ $A^1=\da \{x\in P: \exists $\ directed set\ $D\subseteq \da A,\ \bigvee D=x\}$.

\begin{lemma}\label{lem quasi-c}{\rm \cite{GK2017}}
Let\ $P$\ be a quasi-continuous\ domain\ and\ $A\subseteq P$.\ Then,\ the Scott\ closure of\ $A$\ equals to\ $A^1$,\ that is\ $\overline{A}=A^1$.
\end{lemma}

\begin{corollary}{\rm\cite{XXLK2020}}
Let\ $P$\ be a continuous or quasi-continuous domain.\ Then\ $\sigma(C(X))|_{LP}=O_{\Rightarrow_L}(LP)$,\ that is,\ endowed with the Scott\ topology,\ $P$\ is a directed space,\ the directed lower powerspace\ $P_L(X)$\ is a subspace of the lower powerdomain\ $H(P)$.
\end{corollary}
\noindent{\bf Proof} By Proposition\ \ref{prop rstLX},\ we only need to prove that,\ $\forall \mathcal{V}\in O_{\Rightarrow_L}(LP)$,\ $\uparrow_{C(P)}\mathcal{V}\in\sigma(C(P))$.\ Suppose\ ${\mathcal D}\subseteq C(P)$\ and\ $\overline{\bigcup {\mathcal D}}\in \uparrow_{C(P)}\mathcal{V}$.\ Thus we have some\ $\da F\in \mathcal{V}$\ such that\ $\da F\subseteq \overline{\bigcup {\mathcal D}'}$.\ Let\ ${\mathcal D}'=\{\da F: F\subseteq_{fin}\bigcup{\mathcal D}\}$,\ then\ ${\mathcal D}'$\ is a directed set of\ $LP$\ and\ $\bigcup{\mathcal D}'=\bigcup{\mathcal D}$.\ By Lemma\ \ref{lem quasi-c},\ $\forall a\in F$,\ $\exists$\ a directed set\ $D_a\subseteq \bigcup {\mathcal D}'$\ such that\ $D_a\rightarrow a$.\ Thus,\ ${\mathcal D}'\Rightarrow_L \da F\in \mathcal{V}$.\ Noticed that\ $\mathcal{V}\in O_{\Rightarrow_L}(LP)$,\ so we have some\ $\da G\in {\mathcal D}'\cap \mathcal{V}$.\ For\ ${\mathcal D}$\ is directed,\ there exists\ $A\in {\mathcal D}$\ such that\ $\da G\subseteq A$.\ Thus\ $A\in \uparrow_{C(P)}\mathcal{V}$.\ That is,\ $\uparrow_{C(P)}\mathcal{V}\in\sigma(C(P))$.\ $\Box$

\vskip 3mm

\section{The Directed Upper Powerspaces of Directed Spaces}

In this section,\ we will construct the directed upper powerspace of the directed space,\ which is a free algebra generated by the directed deflationary operation of the directed space.

	\begin{definition}\label{def ddsl}
		Let\ $X$\ be a directed space.
\begin{enumerate}
\item[(1)] A binary operation\ $\oplus :X\otimes X\rightarrow X$\ on\ $X$\ is called a deflationary operation if it is continuous and satisfy the following four conditions:\ $\forall x,\ y,\ z\in X$,\ \begin{enumerate}
		\item[(a)] $x\oplus x=x$,
		\item[(b)] $(x\oplus y)\oplus z=x\oplus (y\oplus z)$,
        \item[(c)] $x\oplus y=y\oplus x$,			
        \item[(d)] $x\oplus y\leq x$.
		\end{enumerate}
\item[(2)] If\ $\oplus$\ is a deflationary operation on\ $X$,\ then\ $(X,\ \oplus)$\ is called a directed deflationary semilattice,\ that is,\ directed deflationary semilattices are those directed spaces with deflationary operations.
\end{enumerate}
	\end{definition}
By Theorem\ \ref{theorem opc}(2),\ the operation\ $\oplus$\ on a directed space\ $X$\ is continuous if
and only if it is monotone and for each given\ $ x,\ y\in X$\ and directed set\ $D\subseteq X$,\ $D\rightarrow x$\ implies\ $(D\oplus y)\rightarrow x\oplus y$.\ Here,\ $D\oplus y=\{d\oplus y:d\in D\}$.

\vskip 3mm
Here are two examples.

\begin{example}
\begin{enumerate}\
 \item[(1)] Suppose\ $P$\ is a poset endowed with the Scott topology,\ and for each\ $a,\ b\in P$,\ the infimum of\ $a$\ and\ $b$\ exists in\ $P$\ (denote by \ $a\wedge b$).\ Then\ $(P,\ \wedge)$\ is a directed deflationary semilattice.
 \item[(2)] Let\ $I=[0,\ 1]$\ (the unit interval),\ let\ ${\mathcal T}$\ denote the topology generated by\ $\{[0,\ a]:a\in I\}$.\ It is easy to check that\ $(I,\ \mathcal T)$\ is a directed space,\
and\ $(I,\ min)$\ is a directed deflationary semilattice endowed with\ $\mathcal T$.
 \end{enumerate}
\end{example}
\vskip 3mm
\begin{definition}
Suppose\ $(X,\ \oplus),\ (Y,\ \uplus)\ $\ are two directed deflationary semilattices,\ function\ $f:(X,\ \oplus)\rightarrow (Y,\ \uplus)$\ is called a deflationary homomorphism
between\ $X$\ and\ $Y$,\ if\ $f$\ is continuous and\ $f(x\oplus y)=f(x)\uplus f(y)$\ holds for\ $\forall x,\ y\in X$.
\end{definition}

Denote the category of all directed deflationary semilattices and deflationary homomorphisms by\ ${\bf Ddsl}$.\ Then ${\bf Ddsl}$\ is a subcategory of\ ${\bf Dtop}$.
\vskip 3mm
\begin{lemma}
Suppose\ $(X,\ \oplus)$\ is a directed deflationary semilattice,\ then\ $\oplus =\wedge_{\sqsubseteq} $.\ Here,\ $\forall x,\ y \in X,\ x\wedge_{\sqsubseteq}y$\ means the infimum of \ $x$\ and\ $y$\ respect to the specialization order\ $\sqsubseteq$\ on\ $X$\ (calling it meet operation).\ Conversely,\ suppose\ $X$\ is a directed space and for each\ $x,\ y \in X $,\ $x\wedge_{\sqsubseteq}y$\ exists,\ the continuity of\ $\wedge_{\sqsubseteq} $\ will naturally imply that\ $(X,\ \wedge_{\sqsubseteq})$\ is a directed deflationary semilattice.

\end{lemma}
\noindent{\bf Proof}  By Definition\ \ref{def ddsl},\ $\forall x,\ y\in X,\ x\oplus y\leq x,\ y$,\ that is\ $x\oplus y$\ is a lower bound
of\ $\{x,\ y\}$.\ Suppose\ $z$\ is another arbitrary lower bound of\ $\{x,\ y\}$,\ by
Theorem\ \ref{theorem opc}(1),\ the pointwise order equals to the specialization order of\ $X\otimes X$,\ then
\ $(x,\ y)\sqsubseteq (z,\ z)$.\ By the continuity and idempotence of the deflationary operation,\
we have\ $z\oplus z=z\sqsubseteq x\oplus y $.\ That is\ $x\oplus y$\ is the infimum of\ $\{x,\ y\}$,\ which
means\ $x\oplus y=x\wedge_{\sqsubseteq}y$.\ Conversely,\ continuous meet operation naturally satisfy
all conditions in Definition\ \ref{def ddsl},\ thus,\ we get the conclusion.\ $\Box$

\vskip 3mm
The above results show that a directed deflationary semilattice\ $(X,\ \oplus)$\ is just a directed space with a continuous meet operation\ $\wedge_{\sqsubseteq}$\ satisfying\ $\oplus=\wedge_{\sqsubseteq}$.\ Since the order on the directed space in this paper is always the specialization order,\ we will use the symbol\ $\wedge$\ instead of\ $\wedge_{\sqsubseteq}$\ in the following part.\ Therefore,\ a directed deflationary semilattice is always represented by a tuple of the form\ $(X,\ \wedge)$,\ here\ $X$\ is a directed space,\ $\wedge$\ represents the continuous meet operation on\ $X$.

\vskip 3mm
Next,\ we give the definition of directed upper powerspace.

\begin{definition}\label{def dups}
Suppose\ $X$\ is a directed space.\ A directed space\ $Z$\ is called the directed upper powerspace over\ $X$\ if and only if the following two conditions are satisfied:

\begin{enumerate}
  \item[(1)] $Z$\ is a directed deflationary semilattice,\ that is the meet operation\ $\wedge$\ on\ $Z$\ exists and which is continuous,
  \item[(2)] There is a continuous function\ $i:\ X\longrightarrow Z$\ satisfy: for an arbitrary directed deflationary semilattice \ $(Y,\ \wedge)$\ and continuous function\ $f:X\longrightarrow Y$,\ there exists
a unique deflationary homomorphism\ $\bar{f}:(Z,\ \wedge)\rightarrow (Y,\ \wedge)$\ such that\ $f=\bar{f}\circ i$.

\end{enumerate}
\end{definition}
By the definition above,\ if directed deflationary semilattices\ $(Z_1,\ \wedge)$\ and\ $(Z_2,\ \wedge)$\ are both the directed upper powerspaces of\ $X$,\ then there exists a topological homomorphism which is also a deflationary homomorphism\ $g:Z_1\rightarrow Z_2$.\ Therefore,\ up to order isomorphism and topological homomorphism,\ the directed upper powerspace of a directed space is unique.\ Particularly,\  we denote the directed upper powerspace of each directed space\ $X$\ by\ $P_U(X)$.

\vskip 3mm
Now,\ we will prove the existence of the directed upper powerspace of each directed space\ $X$\ by way of concrete construction.

\vskip 3mm

	Let\ $X$\ be a directed space.\ Denote
$$UX=\{\uparrow F :F\subseteq _{fin} X\},$$
here,\ $F\subseteq _{fin} X$\ is an arbitrary nonempty finitary subset of\ $X$.\ Define an order\ $\leq_U$\ on\ $UX$\ :
$$\uparrow F_{1}\leq_U \uparrow F_{2}
	\iff \uparrow F_{2}\subseteq \uparrow F_{1}.$$
	Let\ $\mathcal{F}\subseteq UX$\ be a directed set(respect to order\ $\leq_U$)\ and\ $\ua F\in UX$.\ Define a convergence notation\ $\mathcal{F}\Rightarrow_{U}\ua F\iff$\ there exists finite directed sets\ $D_{1},\dots,D_{n}\subseteq X$\ such that
\begin{enumerate}
	\item $F\cap \lim D_{i}\neq \emptyset,\ \forall D_{i}$;
	\item $F\subseteq \bigcup\limits_{i=1}^{n}\lim D_{i}$;
	\item $\forall (d_{1},\dots,d_{n})\in \prod\limits_{i=1}^{n}D_{i}$,\ there exists some\ $\uparrow F^{'}\in \mathcal{F}$,\ such that\ $\uparrow F^{'}\subseteq \bigcup\limits_{i=1}^{n}\uparrow d_{i}$.
\end{enumerate}
A subset\ $\mathcal{U}\subseteq UX$\ is called a\ $\Rightarrow_U$\ convergence open set of\ $UX$\ if and only if for each directed subset\ $\mathcal{F}$\ of\ $UX$\ and\ $\ua F\in UX$,\ $\mathcal{F}\Rightarrow_{U}\uparrow F\in \mathcal{U}$\ implies\ $\mathcal{F}\cap \mathcal{U}\neq \emptyset$.\ Denote all\ $\Rightarrow_U$\ convergence open set of\ $UX$\ by\ $O_{\Rightarrow_{U}}(UX)$.

\begin{proposition}\label{prop UX}
Suppose\ $X$\ is a directed space,\ the following are true:
	\begin{enumerate}
		\item[(1)] $(UX,\ O_{\Rightarrow_U}(UX))$\ is a topological space,\ abbreviated as\ $UX$.
		\item[(2)] The specialization order\ $\sqsubseteq $\ of\ $(UX,\ O_{\Rightarrow_U}(UX))$\ equals to\ $\leq_U$.
		\item[(3)] $(UX,\ O_{\Rightarrow_U}(UX))$\ is a directed space,\ that is\ $O_{\Rightarrow_U}(UX)=d(UX)$.

	\end{enumerate}
\end{proposition}
\noindent{\bf Proof}
	 (1)  Obviously we have\ $\emptyset,\ UX\in O_{\Rightarrow_U}(UX)$.\ If\ $\mathcal{U}\in O_{\Rightarrow_U}(UX)$,\ and\ $\uparrow F_{1}\leq\uparrow F_{2},\uparrow F_{1}\in\mathcal{U},\ F_{1}=\{a_{1},\dots,a_{n}\}$.\ Then it is evident that\ $\{\uparrow F_{2}\}\Rightarrow_{U}\uparrow F_{1}$,\ since we only need to take\ $D_{i}=\{a_{i}\},\ i=1,\dots,n$.\ Then,\ $\{\uparrow F_{2}\}\cap \mathcal{U} \neq
	\emptyset$,\ this means\ $\uparrow F_{2}\in \mathcal{U}$,\ and\ $\mathcal{U}$\ is an upper set respect to order\ $\leq_U$,\ $\mathcal{U}=\uparrow_{\leq}\mathcal{U}$.
	
	Let\ $\mathcal{U}_{1},\mathcal{U}_{2}\in O_{\Rightarrow_U}(UX)$,\ and a directed set \ $\mathcal F\subseteq UX$\ with\ $\mathcal{F}\Rightarrow_{U} \uparrow F\in \mathcal{U}_{1}\cap\mathcal{U}_{2}$,\ then,\ there exists\ $\uparrow F_{1}\in \mathcal{F}\cap \mathcal{U}_{1}$\ and\ $
	\uparrow F_{2}\in \mathcal{F}\cap \mathcal{U}_{2}$,\ but\ $\mathcal{F}$\ is directed,\ we have\ $\uparrow F_{3}\in \mathcal{F},\ \uparrow F_{3}\subset\uparrow F_{1}$\ and\ $\uparrow F_{2}$.\ Then,\ $\uparrow F_{3}\in \mathcal{F}\cap \mathcal{U}_{1}\cap\mathcal{U}_{2}$.\ By the same way,\ we can prove that \ $O_{\Rightarrow_U}(UX)$\ is closed under arbitrary union.\ It follows that\ $O_{\Rightarrow_U}(UX)$\ is a topology.

\vskip 3mm
\
(2) Let\ $\ua F_1,\ \ua F_2\in UX$.\ If\ $\ua F_{1}\leq_L \ua F_{2}$,\ By the proof of\ (1),\ each\ $\Rightarrow_U$\ convergence open set is an upper set respect to\ $\leq_U$,\ then\ $\ua F_1\in \overline{\{\ua F_2\}}$,\ that is,\ $\ua
	        F_{1}\sqsubseteq \ua F_{2}$.

     On the other hand,\ suppose\ $\ua F_{1}\sqsubseteq \ua F_{2}$.\ We need to prove that\ $\ua F_{1}\leq_U \ua F_{2}$,\ that is to prove that\ $\{\ua F\in UX: \ua F_2\subseteq \ua F\}$\ is a closed set in \ $UX$\ respect to topology\ $O_{\Rightarrow_U}(UX)$,\ since\ $\{\ua F\in UX: \ua F_2\subseteq \ua F\}=\{\ua F:\ua F\leq_U\ua F_2\}\subseteq\{\ua F:\ua F\sqsubseteq \ua F_2\}= \overline{\{\ua F_2\}}$,\ then,\ $\{\ua F\in UX: \ua F_2\subseteq \da F\}$\ is a closed set in\ $UX$\ respect to\ $O_{\Rightarrow_U}(UX)$\ implies\ $\{\ua F\in UX: \ua F\subseteq \ua F_2\}=\overline{\{\ua F_2\}}$,\ that is\ $\uparrow F_{1}\leq_U \uparrow F_{2}$.\ Now,\ we prove that\ $\{\ua F\in UX: \ua F_2\subseteq \ua F\}$\ is closed in\ $UX$\ respect to\ $O_{\Rightarrow_U}(UX)$,\ equivalently,\ $\mathcal{U}=UX\setminus\{\ua F\in UX: \ua F_2\subseteq \ua F\}$\ is a\ $\Rightarrow_U$\ convergence open set in\ $UX$.

By contradiction,\ suppose\ $\mathcal{U}$\ is not a\ $\Rightarrow_U$\ convergence open set.\ Then there exists a directed set\ $\mathcal F$\ of\ $UX$\ with\ $\mathcal{F}\Rightarrow_{U} \uparrow F \in \mathcal{U}$\ but\ $\mathcal{U}\cap\mathcal{F}=\emptyset$.\ According to the definition of\ $\Rightarrow_U$\ convergence,\ there exists finite directed sets\ $D_{1},\dots,\ D_{n}\subseteq X$\ such that
\begin{enumerate}
	\item $F\cap \lim D_{i}\neq \emptyset,\ i=1,\ 2,\ldots,n$;
	\item $F\subseteq \bigcup\limits_{i=1}^{n}\lim D_{i}$;
	\item $\forall (d_{1},\dots,d_{n})\in \prod\limits_{i=1}^{n}D_{i}$,\ there exists some\ $\uparrow F^{'}\in \mathcal{F}$,\ such that\ $\uparrow F^{'}\subseteq \bigcup\limits_{i=1}^{n}\uparrow d_{i}$.
\end{enumerate}
Since\ $\uparrow F\in\mathcal{U}$,\ then\ $\uparrow F_{2}\nsubseteq \uparrow F$,\ there exists some\ $a\in F_2$\ with\ $a\not\in\ua F$,\
then\ $F\subseteq X\backslash\da a$.\ According to 1,\ 2 in the definition above,\ for arbitrary\ $i\in\{1,\ 2,\ \ldots,\ n\}$,\ $D_i\cap (X\backslash\da a)\not=\emptyset$.\ For each\ $i$,\ pick\ $d_i\in D_i\cap (X\backslash\da a)$.\ Then\ $(d_1,\ d_2,\ \ldots,d_n)\in \prod\limits_{i=1}^{n}D_{i}$\ and\ $a\not\in \bigcup\limits_{i=1}^n\ua d_i$.\ Since\ ${\mathcal F}\cap {\mathcal U}=\emptyset$,\ we have\ $\forall F'\in{\mathcal F}$,\ $\ua F_2\subseteq \ua F'$,\ thus\ $\ua F'\not\subseteq\bigcup\limits_{i=1}^n\ua d_i$,\ \ this contradicts with\ 3\ in the definition above.\ Therefore,\ $\mathcal{U}$\ is a $\Rightarrow_U$\ convergence open set in\ $UX$.

(3) For an arbitrary topological space\ $X$,\ $O(X)\subseteq d(X)$\ holds,\ then\ $O_{\Rightarrow_U}(UX)\subseteq d(UX)$.\ On the other hand,\ according to the definition of\ $\Rightarrow_U$\ convergence topology,\ if directed set\ $\mathcal F\subseteq UX$\ with\ $
	          \mathcal{F}\Rightarrow_{U}\uparrow F$,\ then\ $\mathcal{F}$\ convergents to\ $\ua F$ respect to\ $O_{\Rightarrow_U}(UX)$.\ Thus,\ by the definition of directed open set,\ $\mathcal{F}\Rightarrow_{F}\uparrow F \in\mathcal{U}\in d(UX)$\ will imply\ $\mathcal{U}\cap \mathcal{F}\neq \emptyset$.\ Then,\ $\mathcal{U}\in O_{\Rightarrow_U}(UX)$,\ it follows that\ $O_{\Rightarrow_U}(UX)= d(UX)$,\ that is,\ $(UX,\   O_{\Rightarrow_U}(UX))$\ is a directed space.\ $\Box$

\begin{proposition}\label{prop UXdc}
 Suppose\ $X,\ Y$\ are two directed spaces.\ Then function\ $f:(UX,\   O_{\Rightarrow_U}(UX))\rightarrow Y$\ is continuous if and only if for each directed set\ ${\mathcal F}\subseteq UX$\ and\ $\ua F\in UX$,\ ${\mathcal F}\Rightarrow_U\ua F$\ implies\ $f(\mathcal{F})\rightarrow f(\ua F)$.
\end{proposition}

\noindent{\bf Proof} Since\ $\Rightarrow_U$\ convergence will lead to\ $O_{\Rightarrow_U}(UX)$\ topological convergence,\ the necessity is obvious.\ We are going to prove the sufficiency.\ Firstly,\ we check
that\ $f$\ is monotone.\ If\ $\ua F_1,\ \ua F_2\in UX$\ and\ $\ua F_1\leq_U\ua F_2$,\ then\ $\{\ua F_2\}\Rightarrow_U\ua F_1$,\ by the hypothesis,\ $\{f(\ua F_2)\}\rightarrow f(\ua F_1)$,\ thus\ $f(\ua F_2)\sqsubseteq f(\ua F_1)$.\ Suppose\ $U$\ is an open set of\ $Y$\ and the directed set\ ${\mathcal F}\Rightarrow_U\ua F\in f^{-1}(U)$,\ then\ $f(\mathcal{F})$\ is a directed set of\ $Y$\ and\ $f(\mathcal{F})\rightarrow f(\ua F)\in U$,\ thus\ $\exists \ua F\in {\mathcal D}$\ such that\ $f(\ua F)\in U$.\ That is,\ $\ua F\in {\mathcal F}\cap f^{-1}(U)$.\ According to the definition of\ $\Rightarrow_U$\ convergence open set,\ $f^{-1}(U)\in O_{\Rightarrow_U}(UX)$,\ that is\ $f$\ is continuous.\ $\Box$
\vskip 3mm

Define a binary operation\ $\cup$\ on\ $UX:\forall \ua F_1,\ \ua F_2\in UX,\ \ua F_1\cup\ua F_2=\ua(F_1\cup F_2)$.\ Suppose we have\ $\ua{F_1}=\ua{F_2},\ \ua{G_1}=\ua{G_2}$\ with\ $F_1\neq F_2,\ G_1\neq G_2$.\ $F_1\cup G_1\subseteq\ua(F_2\cup G_2)$\ implies\ $\ua(F_1\cup G_1)\subseteq\ua(F_2\cup G_2)$.\ Similarly,\ we have the opposite containment,\ thus,\ $\cup$\ is well-defined.
\begin{theorem} Let\ $X$\ be a directed space.\ Then\ $(UX,\ O_{\Rightarrow_U}(UX))$\ respect to the set union operation\ $\cup $\ is a directed deflationary semilattice.
\end{theorem}

\noindent{\bf Proof}
By Proposition\ \ref{prop UX},\ $(UX,\ O_{\Rightarrow_U}(UX))$\ is a directed space.\ We will prove that\ $\cup $\ is a deflationary operation.\ For arbitrary\ $\ua F_1,\ \ua
F_2\in UX$,\ $\ua F_1\cup \ua F_2=\ua (F_1\cup F_2)\in UX$.\ Obviously,\ $\cup$\ satisfy the conditions\ (a),\ (b),\ (c),\ (d)\ in Definition\ \ref{def ddsl},\ we now prove the continuity of\ $\cup$.\
The monotonicity of\ $\cup$\ is evident.\ By Theorem\ \ref{theorem opc}(2)\ and Proposition\ \ref{prop UXdc},\ we only need to prove that,\ for each directed set\ ${\mathcal F}\subseteq UX$\ and\ $\ua F,\ \ua G\in UX$,\ ${\mathcal F}\Rightarrow_U\ua F$\ will imply\ $G\cup {\mathcal F}\Rightarrow_U\ua G\cup \ua F=\ua(G\cup F)$.\ Here,\ $G\cup {\mathcal F}=\{\ua (G\cup F'):\ua F'\in \mathcal{F}\}$\ is still a directed set.\ According to the definition of\ $\Rightarrow_U$\ convergence,\ there exists finite directed sets\ $D_{1},\dots,D_{k}\subseteq X$,\ satisfy the conditions such that\ ${\mathcal F}\Rightarrow_U\ua F$.\ Let\ $G=\{a_{1},\dots,a_{n}\}$,\ and\ $D_{k+1}=\{a_{1}\},\ D_{k+2}=\{a_2\},\ \ldots,\ D_{k+n}=\{a_{n}\}$.\ It is straightward to verify that,\ $D_1,\ D_2,\ \ldots,\ D_k,\ D_{k+1},\ \ldots,\ D_{k+n}$\ satisfy all the conditions such that\ $G\cup {\mathcal F}\Rightarrow_U(\ua G\cup \ua F)$.\ It follows that,\ $(UX,\ \cup)$\ is a directed deflationary semilattice.\ $\Box$

\vskip 3mm
The following theorem is the main result in this paper.

\begin{theorem}\label{theorem UX}
	Suppose\ $X$\ is a directed space,\ then\ $(UX,\   O_{\Rightarrow_U}(UX))$\ is the directed upper powerspace over\ $X$,\ that is,\ endowed with topology\ $O_{\Rightarrow_U}(UX)$,\ $(UX,\ \cup)\cong P_U(X)$.	
\end{theorem}

{\bf Proof} Define function\ $i:X\rightarrow UX$:\ $\forall x\in X$,\ $i(x)=\uparrow x$.\ We prove the continuity of\ $i$.\ It is evident that\ $i$\ is monotone.\ Suppose we have directed set\ $D\subseteq X$\ and\ $x\in X$\ satisfy\ $D\rightarrow x$.\ Let\ $\mathcal{D}=\{\ua d:d\in D\}$,\ then\ $\mathcal{D}$\ is a directed set in\ $UX$\ and\ ${\mathcal D}\Rightarrow_{U}\uparrow x$.\ Notice that\ $i(D)=\mathcal{D}$,\ so\ $i(D)\Rightarrow_{U}\ua x=i(x)$.\ By Proposition\ \ref{prop dc},\ $i$\ is continuous.

Let\ $(Y,\ \wedge)$\ be an arbitrary directed deflationary semilattice,\ $f:X\rightarrow Y$\ is a continuous function.\ Define\ $\bar{f}:UX\rightarrow Y$\ as follows:\ $\forall \ua F\in UX$ (let\ $F=\{a_1,\ a_2,\ \ldots,\ a_n\}$),\ $$\bar{f}(\ua F)=f(a_1)\wedge f(a_2)\wedge\cdots\wedge f(a_n)=\bigwedge\limits_{a\in F}f(a).$$
Particularly,\ we write\ $\bar{f}(\ua F)=\wedge f(F)$.\ $\bar{f}$\ is well-defined,\ since if we have\ $\ua F =\ua G$\ with\ $F\neq G$,\ then\ $f(F)\subseteq \ua f(G)$\ implies\ $\wedge f(G)\leq\wedge f(F)$,\ that is\ $\bar{f}(\ua G)\leq\bar{f}(\ua F)$.\ Similarly,\ we have\ $\bar{f}(\ua F)\leq\bar{f}(\ua G)$.\ Thus,\ $f$ is well-defined.

(1) $f=\bar{f}\circ i$.

For arbitrary\ $x\in X$,\ $(\bar{f}\circ i)(x)=\bar{f}(i(x))=\bar{f}(\ua x)=f(x)$.

(2) $\bar{f}$\ is a deflationary homomorphism,\ that is,\ $\bar{f}$\ is continuous and for arbitrary\ $\ua F_1,\   \ua F_2\in UX$,\ $\bar{f}(\ua F_1\cup \ua F_2)=\bar{f}(\ua F_1)\wedge \bar{f}(\ua F_2)$.

First,\ we prove that\ $\bar{f}$\ preserves the union operation.\ Suppose\ $\ua F_1,\ \ua F_2\in UX$.\ Then\ $\bar{f}(\ua F_1\cup\ua F_2)=\bar{f}(\ua(F_1\cup F_2))=\wedge f(F_1\cup F_2)=(\wedge f(F_1)\wedge(\wedge f(F_2))=\bar{f}(\ua F_1)\wedge\bar{f}(\ua F_2)$.\ Next,\ we prove the continuity of\ $\bar{f}$.\ Notice that\ $\wedge$\ is the meet operation,\ $\bar{f}$\ is evidently monotone.\ Suppose\ ${\mathcal F}\subseteq UX$\ is a directed set and\ $\ua F\in UX$\ satisfy\ ${\mathcal F}\Rightarrow_U\ua F$.\ By the definition of\ $\Rightarrow_U$,\ there exists finite directed sets\ $D_{1},\dots,D_{n}\subseteq X$\ such that
\begin{enumerate}
	\item $F\cap \lim D_{i}\neq \emptyset,\ i=1,2,\ldots,n$;
	\item $F\subseteq \bigcup\limits_{i=1}^{n}\lim D_{i}$;
	\item $\forall (d_{1},\dots,d_{n})\in \prod\limits_{i=1}^{n}D_{i}$,\ there exists some\ $\uparrow F^{'}\in \mathcal{F}$,\ such that\ $\uparrow F^{'}\subseteq \bigcup\limits_{i=1}^{n}\uparrow d_{i}$.
\end{enumerate}

Let\ $F=\{b_1,\ b_2,\ldots,b_k\}$.\ By\ 1,\ for each\ $1\leq i\leq n$,\ we have some\ $b_i\in F$\ such that\ $D_i\rightarrow b_i$.\ If\ $F\setminus \{b_1,\dots,b_n\}\neq\emptyset$,\ which is denoted by\ $G=\{a_1,\   a_2\ldots,\ a_s\}$.\ By\ 2,\ For each\ $a_j\in G$,\ we have\ $1\leq i_j\leq n$\ such  that\ $D_{i_j}\rightarrow a_j$.\ By the continuity of\ $f$,\ then\ $f(D_{i})\rightarrow f(b_{i}),\ i=1,\dots,n$\ and,\ $f(D_{i_j})\rightarrow f(a_j)$,\ $j=1,2,\ldots,s$.\ Since the meet operation\ $\wedge$\ on\ $Y$\ is continuous,\ the following convergence holds
$$f(D_{1})\wedge\cdots \wedge f(D_n)\wedge f(D_{i_1})\wedge\cdots\wedge f(D_{i_s})\rightarrow f(b_{1})\wedge \dots \wedge f(b_n)\wedge f(a_{i_1})\wedge\cdots\wedge f(a_{i_s}).\ \\ (\ast)$$
Here,\ $f(D_{1})\wedge \cdots \wedge f(D_n)\wedge  f(D_{i_1})\wedge\cdots\wedge f(D_{i_s})=\{f(d_{1}) \wedge \cdots \wedge f(d_n)\wedge f(d_{i_1})\wedge\cdots\wedge f(d_{i_s}) :(d_1,\   \ldots,\   d_k,\   d_{i_1},\ \ldots,\   d_{i_s})\in (\prod\limits_{i=1}^{n}D_i\})\times(\prod\limits_{j=1}^{s}D_{i_j})\}$.\ Let\ $U$\ be an arbitrary open neighborhood of\ $\wedge f(F)$,\ by\ $(\ast)$,\ there exists\ $(d_1,\ \ldots,\ d_n,\ d_{i_1},\ \ldots,\ d_{i_s})\in (\prod\limits_{i=1}^{n}D_i\})\times(\prod\limits_{j=1}^{s}D_{i_j})$\ such that\ $f(d_{1})\wedge
 \cdots \wedge f(d_{n})\wedge f(d_{i_1})\wedge\cdots\wedge f(d_{i_s})\in U$.\ For each\ $D_{i_j}$\ repeats\ $D_i$,\ and each\ $D_i$\ is directed,\ we have\ $(d^{'}_1,\ d^{'}_2,\ \ldots,\ d^{'}_n)\in\prod\limits_{i=1}^{n}D_i$\ such that\ $f(d^{'}_1)\wedge\cdots\wedge f(d^{'}_n)\sqsupseteq f(d_{1})\wedge
 \cdots \wedge f(d_{n})\wedge f(d_{i_1})\wedge\cdots\wedge f(d_{i_s})$.\ By\ 3,\ there exists some\ $\uparrow F^{'}\in\mathcal{F}$\ such that\ $\uparrow F^{'}\subseteq
 \bigcup\limits_{i=1}^{n}\uparrow d^{'}_{i}$.\ Thus\ $\bar{f}(\uparrow F^{'})=\wedge f(F^{'})\sqsupseteq f(d^{'}_{1})\wedge\dots\wedge f(d^{'}_{n})$.\ But\ $U$\ is an upper set,\ it follows that\ $
 \wedge f(F^{'})=\bar{f}(\uparrow F^{'})\in U$,\ then\ $\bar{f}({\mathcal F})=\{\wedge f(F^{'}):\uparrow F^{'}\in\mathcal{F}\}\rightarrow\wedge f(F)$.\ By Proposition\ \ref{prop UXdc},\ function\ $\bar{f}$\ is continuous.

(3) Homomorphism\ $\bar{f}$\ is unique.

Suppose we have a deflationary homomorphism\ $g:(UX,\ \cup)\rightarrow (Y,\ \wedge)$\ such that\ $f=g\circ i$,\ then\ $g(\uparrow x)=f(x)=\bar{f}(\ua x)$.\ For each\ $ \ua F\in UX$ (let\ $ F=(a_{1},\dots,a_{n})$),\ \begin{eqnarray*}
 g(\uparrow F)&=& g(\ua a_1\cup\ua a_2\cdots\cup\ua a_n)\\ &=&g(\uparrow a_{1})\wedge g(\ua a_2)\wedge \cdots \wedge g(\uparrow a_{n})\\ &=&\bar{f}(\uparrow a_{1})\wedge \bar{f}(\ua a_2)\wedge\cdots\wedge\bar{f}(\uparrow a_{n})\\ &=&\bar{f}(\ua a_1\cup\ua a_2\cdots\cup\ua a_n)\\ &=&\bar{f}(\uparrow F).
 \end{eqnarray*}
 Thus\ $\bar{f}$\ is unique.

 In conclusion,\ according to definition\ \ref{def ddsl},\ endowed with topology\ $O_{\Rightarrow_U}(UX)$,\ the directed deflationary semilattice\ $(UX,\   \cup)$\ is the directed upper powerspace of\ $X$,\ that is,\ $P_U(X)\cong (UX,\   \cup)$.\ $\Box$

 \vskip 3mm
The directed upper powerspace is unique in the sense of order isomorphism and topological homomorphism,\  so we can directly denote the directed upper powerspace by\ $P_U(X)=(UX,\ \cup)$\  of each directed space\ $X$.

 Suppose\ $X,\ Y$\ are two directed spaces,\ $f:X\rightarrow Y$\ is a continuous function.\ Define function\ $P_U(f):P_U(X)\rightarrow P_U(Y)$\ as follows:\ $\forall \ua F\in UX$,\ $$P_U(f)(\ua F)=\ua f(F).$$
 We can check that,\ $P_U(f)$\ is well-defined and order preserving.\ It is easy to check that,\ $P_U(f)$\ is a deflationary homomorphism between these two lower powerspaces.\ If\ $id_X$\ is the identity function
and\ $g:Y\rightarrow Z$\ is an arbitrary continuous function from\ $Y$\ to a directed space\ $Z$,\ then,\ $P_U(id_X)=id_{P_U(X)},\ P_U(g\circ f)=P_U(g)\circ P_U(f)$.\ Thus,\ $P_U:{\bf Dtop}\rightarrow {\bf Ddsl}$\ is a functor from \ ${\bf Dtop}$\ to\ ${\bf Ddsl}$.\ Let\ $U:{\bf Ddsl}\rightarrow {\bf Dtop}$\ be the forgetful functor.\ By Theorem\ \ref{theorem UX},\ we have the following result.
\vskip 3mm
\begin{proposition}\label{P_U morphism}
	$P_U$\ is a deflationary homomorphism between\ $P_U(X)$\ and\ $P_U(Y)$\ for each directed spaces\ $X$\ and\ $Y$.
	
\end{proposition}
{\bf Proof} It is directed to check that\ $P_U$\ preserves deflationary operation between\ $P_U(X)$\ and\ $P_U(Y)$\ for each directed spaces\ $X$\ and\ $Y$.\ We only to chek that\ $P_U$\ is continuous.\

Let\ $X$\ and\ $Y$\ be two directed spaces,\ let\ $\mathcal{D}=\{f_i\}_{i\in I}$\ be a directed set in\ $[X\rightarrow Y]$\ and convergent to\ $f\in[X\rightarrow Y]$\ pointwise.\ Let\ $\uparrow F\in P_U(X)$\ and\ $F=\{a_1,\dots,a_n\}$,\ we shall verify that\ $\{P_U(f_i)(\uparrow F)\}_{i\in I}\rightarrow P_U(f)(\uparrow F)$,\ that is\ $\{\uparrow f_i(F)\}_{i\in I}\rightarrow\uparrow f(F)$\ in\ $P_U(Y)$.\ According to the definition of\ $\Rightarrow_{U}$\ convergence,\ firstly,\ we have\ $D_a=\{f_i(a)\}_{i\in I}\rightarrow f(a)$\ in\ $Y$\ for each\ $a\in F$.\ Secondly,\ for each\ $f_{i_1}(a_1)\in D_{a_1},\dots,f_{i_n}(a_n)\in D_{a_n}$,\ since each\ $a\in F,\ D_a$\ is directed,\ we have some\ $\uparrow f_{i_0}(F)\in\{\uparrow f_i(F)\}_{i\in I}$\ such that\ $\uparrow f_{i_0}(F)\subseteq\uparrow(f_{i_1}(a_1),\dots,f_{i_n}(a_n))$.\ Thus\ $\{\uparrow f_i(F)\}_{i\in I}\Rightarrow\uparrow_U f(F)$,\ $P_U$\ is continuous.

 \begin{corollary}
 $P_U$\ is a left adjoint of the forgetful functor\ $U$,\ that is,\ ${\bf Ddsl}$\ is a reflective subcategory of\ ${\bf Dtop}$.
 \end{corollary}

\section{Relations Between Upper Powerspace}

In this section,\ we will discuss the relations between the upper powerdomains of dcpos,\ the observationally-induced upper powerspaces and the directed upper powerspaces.

According to the results in the last section,\ for an arbitrary directed space,\ the directed lowerspace is the set of\ $UX$,\ endowed with the\ $\Rightarrow_U$\ convergence topology,\ the deflationary operation equals to the union operation of sets.\ In general,\ for an arbitrary topological space and arbitrary\ dcpo,\ although the observationally-induced upper powerspace and the upper powerdomain exist,\  their concrete structure cannot be expressed(see article\ \cite{BS2015,  Hec91,  Hec92}).

Let\ $(X,\ O(X))$\ be a topological space.\ We say that an nonempty set\ $A\subseteq X$\ is a saturated set,\ if\ $A=\bigcap\{U\in O(X):A\subseteq U\}$.\ Denote all nonempty compact saturated sets of\ $X$\ by\ $Q(X)$.\ For each\ $U\in O(X)$,\ let
$$[U]=\{K\in Q(X): K\subseteq U\}.$$
Denote ${\mathcal B}_X=\{[U]:U\in O(X)\}$.\ The upper\ Vietoris\ topology of\ $Q(X)$\ is generated by subbase\ ${\mathcal {B}_X}$,\ denote by\ $V_U(Q(X))$.\ Particularly,\ for a\ dcpo,\ endowed with the\ Scott\ topology,\ the compact saturated sets are just the compact upper sets.

\begin{theorem}{\rm\cite{BS2015}}\label{theorem oiu}
Suppose\ $X$\ is a\ sober\ and locally compact space,\ $(Q(X),\ V_U(Q(X))$\ is order isomorphic and topological homomorphic to the observationally-induced upper space of\ $X$\ respect to the union operation of sets.\ Under this condition,\ we have\ $V_U(Q(X))=\sigma(Q(X)$.\ Here,\ $Q(X)$\ with the order reverse to containment,\ $\sigma(Q(X)$\ denotes the\ Scott\ topology.
\end{theorem}

\begin{theorem}{\rm\cite{GHK}}
Let\ $P$\ be a continuous domain.\ Then\ $(Q(P),\ \supseteq)$,\ endowed with the\ Scott\ topology and the order reverse to containment,\ isomorphic to the upper powerdomain over\ $P$\ (which is also called\ Smyth\ powerdomain\ )\ $S(P)$,\ besides,\ the following holds:
\begin{enumerate}
\item[(1)] $(Q(P),\ \supseteq)$\ is a continuous meet semilattice,\
 \item[(2)] $\forall K\in Q(P)$,\ $K=\bigcap\{\ua F: 1\leq |F|<\omega\ \&\ K\subseteq (\ua F)^{\circ}\}$.
\end{enumerate}
\end{theorem}

 According to these two theorems above,\ for each continuous domain with the Scott topology,\ its observationally-induced upper powerspace isomorphic to the upper powerdomain.\ Next,\ we discuss the directed upper powerspace of continuous domain.

 Let\ $X$\ be a continuous\ domain.\ Then\ $(X,\ \sigma(X))$\ is a directed space,\ and each upper set of nonempty finitary elements is a compact saturated set.\ Thus,\ $UX\subseteq Q(X)$.\ By\ $\sigma(Q(X))|_{UX}$\ denote the relative topology from the Scott topology on\ $Q(X)$.
 \begin{proposition}\label{prop rstUX}
 Suppose\ $X$\ is a continuous domain which is endowed with Scott topology\ $\sigma(X)$.
 \begin{enumerate}
 \item[(1)] For each given directed set\ ${\mathcal F}\subseteq UX$\ and\ $\ua F\in UX$,\ ${\mathcal F}\Rightarrow_U\ua F\ \Leftrightarrow\ \bigcap\{\ua G:\ua G\in {\mathcal F}\}\subseteq \ua  F$.
 \item[(2)] $O_{\Rightarrow_U}(UX)=\sigma(Q(X))|_{UX}$.
 \end{enumerate}
 \end{proposition}

 \noindent{\bf Proof} (1) Suppose\ ${\mathcal F}\subseteq UX$\ is a directed set of\ $UX$,\ $\ua F\in UX$,\ and\ ${\mathcal F}\Rightarrow_U\ua F$.\ By definition,\ there exist finitary directed sets\ $D_{1}\dots,D_{n}\subseteq X$\ such that
\begin{enumerate}
	\item $F\cap \lim D_{i}\neq \emptyset,\ i=1,2,\ldots,\ n$;
	\item $F\subseteq \bigcup\limits_{i=1}^{n}\lim D_{i}$;
	\item $\forall (d_{1},\dots,d_{n})\in \prod\limits_{i=1}^{n}D_{i}$,\ there exists some\ $\uparrow F^{'}\in \mathcal{F}$,\ such that\ $\uparrow F^{'}\subseteq \bigcup\limits_{i=1}^{n}\uparrow d_{i}$.
\end{enumerate}
By contradiction,\ suppose\ $\bigcap\{\ua G:\ua G\in {\mathcal F}\}\not\subseteq \ua  F$.\ Then there exists some\ $a\in\bigcap\{\ua G:\ua G\in {\mathcal F}\}$\ such that\ $a\not\in\ua F$.\ Thus,\ $F\subseteq X\backslash\da a$.\ By\ 1\ and\ 2,\ for each\ $i$,\ we have\ $d_i\in D_i$\ such that\ $d_i\in X\backslash \da a$.\ By\ 3,\ there exists some\ $\ua F^{'}\in{\mathcal F}$\ such that\ $\ua F^{'}\subseteq \bigcup\limits_{i=1}^n\ua d_i$,\ this contradicts with\ $a\not\in\bigcup\limits_{i=1}^n\ua d_i$.\ That is\ $\bigcap\{\ua G:\ua G\in {\mathcal F}\}\subseteq \ua  F$.

On the other hand,\ suppose\ $\bigcap\{\ua G:\ua G\in {\mathcal F}\}\subseteq \ua  F$.\ Let\ $F=\{a_1,\ a_2,\ldots,a_n\}$.\ For\ $X$\ is a continuous domain,\ then each\ $\dda a_i$\ is directed and\ $a_i=\bigvee \dda a_i$,\ $i=1,2,\ldots,n$.\ Let\ $D_i=\dda a_i$,\ $D_i\rightarrow a_i$,\ then,\ each\ $D_i$\ satisfy\ 1\ and\ 2\ in the definition of\ $\Rightarrow_U$.\ For arbitrary\ $(d_1,d_2,\ldots,\ d_n)\in\prod\limits_{i=1}^nD_i$,\ we have\ $\ua F\subseteq \bigcup\limits_{i=1}^n \dua d_i=(\bigcup\limits_{i=1}^n \ua d_i)^{\circ}$.\ Since each continuous domain\ is\  well-filtered,\ it follows that there exists some\ $\ua G\in{\mathcal F}$\ such that\ $\ua G\subseteq (\bigcup\limits_{i=1}^n \ua d_i)^{\circ}$,\ then,\ 3\ in the definition of\ $\Rightarrow_U$\ holds.\ In conclusion,\ ${\mathcal F}\Rightarrow_U\ua F$.

\noindent (2) Suppose\ ${\mathcal U}\in O_{\Rightarrow_U}(UX)$.\ Let\ ${\mathcal U}_Q=\{K\in Q(X):\exists \ua F\in {\mathcal U},\ K\subseteq\ua F\}$.\ Obviously,\ ${\mathcal U}= {\mathcal U}_Q\cap Q(X)$.\ Let\ ${\mathcal K}\subseteq Q(X)$\ be a directed set respect to the order reverse to containment and\ $\bigcap\{K:K\in{\mathcal K}\}\in{\mathcal U}_Q$.\ There exists some\ $\ua F\in {\mathcal U}$\ such that\ $\bigcap\{K:K\in{\mathcal K}\}\subseteq \ua F$.\ Let\ $F=\{a_1,a_2,\ldots,\ a_n\}$.\ Since\ $X$\ is a continuous domain,\ each\ $\dda a_i$\ is directed and\ $a_i=\bigvee \dda a_i$,\ $i=1,2,\ldots,n$.\ Let\ ${\mathcal F}=\{\bigcup_{i=1}^n\ua d_i:d_i\ll a_i,\ i=1,2,\ldots,n\}$.\ For\ $X$\ is a continuous domain,\ then\ $Q(X)$\ is a continuous domain,\ we have\ $\ua F\subseteq (\bigcup_{i=1}^n\ua d_i)^{\circ}\subseteq \bigcup_{i=1}^n\ua d_i$,\ $i=1,2,\ldots,\ n$,\ besides,\ $\cap\mathcal F=\ua F$.\ By\ (1),\ ${\mathcal F}\Rightarrow_U \ua F$.\ Thus,\ we have\ $(d_1,\   d_2,\   \ldots d_n)\in\prod_{1\leq i\leq n}\dda a_i$\ such that\ $\bigcup_{i=1}^n\ua d_i\in {\mathcal U}$.\ Notice that\ $\bigcap\{K:K\in{\mathcal K}\}\subseteq \ua F\subseteq (\bigcup_{i=1}^n\ua d_i)^{\circ}$,\ there exists some\ $K\in{\mathcal K}$\ such that\ $K\subseteq  (\bigcup_{i=1}^n\ua d_i)^{\circ}$.\ By the definition of\  ${\mathcal U}_Q$,\ $K\in{\mathcal U}$,\ that is,\ ${\mathcal U}_Q$\ is a Scott open set in\ $Q(X)$.\ Therefore,\ $O_{\Rightarrow_U}(UX) \subseteq \sigma(Q(X))|_{UX}$.

On the other hand,\ let\ ${\mathcal V}\in\sigma(Q(X))$,\ and\ ${\mathcal F}\subseteq UX$\ is a directed set with \ ${\mathcal F}\Rightarrow_U\ua F\in {\mathcal V}\cap UX$.\ By\ (1),\ $\bigcap\{\ua G:\ua G\in{\mathcal F}\}\in {\mathcal V}$,\ there exists some\ $\ua G\in{\mathcal F}\cap{\mathcal V}$,\ that is,\ ${\mathcal F}\cap {\mathcal V}\cap UX\not=\emptyset$.\ Thus,\ $ \sigma(Q(X))|_{UX}\subseteq O_{\Rightarrow_U}(UX)$.\ We get the conclusion.\ $\Box$

\begin{example}\label{example UXneq}
Let\ $X={\mathbb R}^n$\ be the\ $n-$dimensional Euclidean space.\ Then\ $X$\ is locally compact\ $T_2$\ space,\ naturally\ sober\ and locally compact.\ Denote the observationally-induced upper space over\ $X$\ by\ $P_O(X)$.\ By Theorem\ \ref{theorem oiu},\ $P_O(X)=\{K\subseteq X: K$\ is a nonempty compact set of\ $X$\ $\}$\ and the topology of which is the Scott topology.\ It is easy to check that,\ $P_O(X)$\ is a continuous domain and for each nonempty compact set \ $K\subseteq X$,\ $\{\{a\}:a\in K\}$\ is a compact saturated set of\ $P_O(X)$.\ It follows that,\ for directed space\ $P_O(X)$,\ its directed upper powerspace\ $P_U(P_O(X))\not=Q(P_O(X))$.

\end{example}

\vskip 3mm

\section{The Directed Convex Powerspace}
In this section,\ we will construct the directed convex powerspace of the directed space,\ which is a free algebra generated by the directed semilattice operation of the directed space.

\begin{definition}\label{def dsl}
		Let\ $X$\ be a directed space.
\begin{enumerate}
\item[(1)] A binary operation\ $ + :X\otimes X\rightarrow X$\ on\ $X$\ is called a semilattice operation if it is continuous and satisfy the following three conditions:: $\forall x,\ y,\ z\in X$,\
		\begin{enumerate}
\item[(a)] $x + x=x$,\ 			
\item[(b)] $(x+ y)+z=x+(y+ z)$,
\item[(c)] $x+y=y +x$.
  		\end{enumerate}
\item[(2)] If\ $ +$\ is a directed semilattice operation on\ $X$,\ then\ $(X,\ +)$\ is called a directed semilattice,\ that is,\ directed  semilattices are those directed spaces with semilattice operations.
\end{enumerate}
	\end{definition}

By Theorem\ \ref{theorem opc}(2),\ the operation\ $+$\ on\ $X$\ is continuous if and only if it is monotone and for each given\ $ x,\ y\in X$\ and directed set\ $D\subseteq X$,\ $ D\rightarrow x$\ implies\ $(D+ y)\rightarrow x + y$.\ Here,\ $D+ y=\{d+y:d\in D\}$.

\vskip 3mm

\begin{example}\rm\ \begin{enumerate}
 \item[(1)] Suppose\ $(P,\ +)$\ is a directed inflationary(or deflationary)semilattice,\ then,\ $(P,\ +)$\ is a directed semilattice.
 \item[(2)] Let\ $X=\{0,\ 1,\ 2\}$\ and\ $0\leq 1\leq 2$,\ define\ $0+0=0,\ 1+1=1,\ 2+2=2,\ 0+1=1,\ 1+2=1,\ 0+2=1$.\ Then\ $(X,\ +\ )$,\ endowed with the Scott topology,\ is a directed semilattice.

 \end{enumerate}
 \end{example}

\begin{definition}
Suppose\ $(X,\ \oplus),\ (Y,\ \uplus)$\ are two directed semilattice,\  function\ $f:(X,\ \oplus)\rightarrow (Y,\ \uplus)$\ is called a directed semilattice homomorphism between\ $X$\ and\ $Y$,\ if\ $f$\ is continuous respect to these two directed space and\ $f(x\oplus y)=f(x)\uplus f(y)$\ holds,\ $\forall x,\ y\in X$.
\end{definition}

\vskip 3mm
Denote the category of all directed semilattice and directed homomorphism by\ ${\bf Dsl}$.\ Then\ ${\bf Dsl}$\ is a subcategory of the category of\ ${\bf Dtop}$.

\vskip 3mm
Next,\ we give the definition of directed convex powerspace.

\begin{definition}\label{def dcps}
  Suppose\ $X$\ is a directed space.\ A directed space\ $Z$\ is called the directed convex powerspace over\ $X$\ if and only if the following two conditions are satisfied:
\begin{enumerate}
\item[(1)] $Z$\ is a directed semilattice,\ that is the semilattice operation\ $+ $\ on\ $Z$\ exists and which is continuous,\
 \item[(2)] There is a continuous function\ $i:X\longrightarrow Z$\ satisfy: for an arbitrary directed  semilattice\ $(Y,\ +)$\ and continuous function\ $f:X\longrightarrow Y$,\ there exists an unique directed semilattice homomorphism\ $\bar{f}:(Z,\ +)\rightarrow (Y,\ +)$\ such that\ $f=\bar{f}\circ i$.
\end{enumerate}
\end{definition}

By the definition above,\ if directed semilattices $(Z_1,\ +)$\ and\ $(Z_2,\ +)$\ are both the directed convex powerspaces of\ $X$,\ then,\ there exists a topological homomorphism which is still a directed semilattice homomorphism\ $g:Z_1\rightarrow Z_2$.\ Therefore,\ up to order isomorphism and topological homomorphism,\  the directed convex powerspace of a directed space is unique.\ Particularly,\ we denote the
directed convex powerspace of each directed space\ $X$\ by\ $P_P(X)$.
\vskip 3mm

Now,\  we will prove the existence of the directed convex powerspace of each directed space\ $X$\ by way of concrete construction.
\vskip 3mm

	Let\ $X$\ be a directed space,\ for each compact set\ $A\subseteq X$,\ $\bar A$\ denotes the closure of\ $A$,\ and\ $sat(A)$\ the saturation of\ $A$,\ that is\ $sat(A)=\cap\{U\in O(X):A\subseteq U\}$.\ Thus,\ for each finite subset\ $F$,\ it follows that\ $\bar F=\da F,\ sat(F)=\ua F$\ respect to the specialization order of\ $X$.\ Denote
$$\widehat{A}=(\bar{A},\  sat(A)),\ \\
PX=\{\widehat{F} :F\subseteq _{fin} X\},\   $$
Here,\ $F\subseteq _{fin} X$\ is an arbitrary nonempty finitary subset of\ $X$.\ Define an order\ $\leq_P$\ on\ $PX$\ as follows:
$$\widehat{F_{1}}\leq_P \widehat{F_{2}}
	\iff \da F_{1}\subseteq \da F_{2}\ \&
 \uparrow F_{2}\subseteq \uparrow F_{1}.$$
	Let\ $\mathcal{D}=\{\widehat{F_{i}}\}_{i\in I}\subseteq PX$\ be a directed set\ (respect to the order\ $\leq_P$),\ and\ $\widehat{F} \in PX$.\ Define a convergence notation\ $\mathcal{D}\Rightarrow_P \widehat{F}$\ as follows:

		$\mathcal{D}\Rightarrow _{P}\widehat{F}  \iff $\ there exists finite directed sets\ $D_{1},\dots,D_{k}$\ of\ $X$\ satisfy the following four conditions:
\begin{enumerate}
\item[(1)] $D_i\subseteq\bigcup\limits_{i\in I}\da F_{i},\ i=1,\dots,k$;
\item[(2)] $\forall i=1,\dots,k,\ F\cap \lim D_{i}\neq \emptyset$;
\item[(3)] $F\subseteq \bigcup\limits_{i=1}^{k}\lim D_{i}$;

\item[(4)] $\forall (d_{1},\dots,d_{k})\in\prod\limits_{i=1}^{k}D_{i},\ \exists\  \widehat{F^{'}}\in\mathcal{D}$,\ such that\ $\ua F^{'}\subseteq \ua(d_{1},\dots,d_{k})$.
\end{enumerate}
A subset\ $\mathcal{U}\subseteq PX$\ is called a\ $\Rightarrow_P$\ convergence open set of\ $PX$\ if and only if for each directed subset\ $\mathcal{D}=\{\widehat{F_{i}}\}_{i\in I}$\ of\ $PX$\ and\ $\widehat{F}\in PX$,\ $\mathcal{D}\Rightarrow_{P}\widehat{F}\in \mathcal{U}$\ will imply\ $\mathcal{D}\cap \mathcal{U}\neq \emptyset$.\ Denote all\ $\Rightarrow_P$\ convergence open set of\ $PX$\ by\ $O_{\Rightarrow_{P}}(PX)$.

\begin{proposition}\label{prop PX}
	Suppose\ $X$\ is a directed space,\ the following are true:
	\begin{enumerate}
		\item[(1)] $(PX,\ O_{\Rightarrow_P}(PX))$\ is a topological space,\ abbreviated as\ $PX$.
		\item[(2)] The specialization order\ $\sqsubseteq $\ on\ $(PX,\ O_{\Rightarrow_P}(PX))$\ equals to\ $\leq_P$.
		\item[(3)] $(PX,\ O_{\Rightarrow_P}(PX))$\ is a directed space,\ that is\ $O_{\Rightarrow_P}(PX)=d(PX)$.
		
	\end{enumerate}
\end{proposition}
\noindent{\bf Proof}
	 (1)  Firstly,\ $\emptyset,\ PX\in O_{\Rightarrow_P}(PX)$.\ Suppose\ $\mathcal{U}\in O_{\Rightarrow_P}(PX)$,\ $\widehat{F_{1}},\ \widehat{F_{2}}\in PX$,\ Let\ $F_{1}=\{a_{1},\dots,a_{n}\}$.\ If\ $\widehat{F_{1}}\leq_{P} \widehat{F_{2}}$,\ then\ $\{\widehat{F_{2}}\}\Rightarrow_{P}\widehat{F_{1}}$\ holds,\ since the finite directed sets\ $\{a_{1}\},\dots,\{a_{n}\}$\ is sufficient to satisfy the definition of \ $\{\widehat{F_{2}}\}\Rightarrow_{P}\widehat{F_{1}}$,\ then\ $\widehat{F_{1}}\in\mathcal{U}$\ implies\ $\widehat{F_{2}}\in\mathcal{U}$.\ It follows that\ $\mathcal{U}$\ is an upper set of\ $PX$\ respect to\ $\leq_P$.\ Suppose\ $\mathcal{U}_{1},\ \mathcal{U}_{2}\in O_{\Rightarrow_P}(PX)$,\ $\mathcal{D}$\ is a directed set in\ $PX$,\ $\widehat{F}\in PX$\ and satisfy\ $\mathcal{D} \Rightarrow_{P} \widehat{F}\in \mathcal{U}_{1}\cap\mathcal{U}_{2}$,\ then there exists some\ $\widehat{G_{1}}\in \mathcal{D}\cap \mathcal{U}_{1}$\ and\ $\widehat{G_{2}}\in \mathcal{D}\cap \mathcal{U}_{2}$.\ Since\ $\mathcal{D}$\ is directed,\ we have\ $\widehat{G}\in \mathcal{D}$\ such that\ $\widehat{G_{1}},\   \widehat{G_{2}}\leq_P \widehat{G} $,\ then\ $\widehat{G}\in \mathcal{D}\cap \mathcal{U}_{1}\cap\mathcal{U}_{2}$.\ By the same way,\ we can evidently prove that\ $O_{\Rightarrow_P}(PX)$\ is closed under arbitrary union.\ That is,\ $O_{\Rightarrow_P}(PX)$\ is a topology on\ $PX$.

\vskip 3mm

(2) Let\ $\widehat{F_{1}},\widehat{F_{2}}\in PX$.\ If\ $ \widehat{F_{1}}\leq_P \widehat{F_{2}}$,\ according the proof of\ (1),\ each\ $\Rightarrow_P$\ convergence open set is an upper set respect to partial order\ $\leq_P$,\ then\ $\widehat{F_{1}}\in \overline{\{\widehat{F_{2}}\}}$,\ that is,\ $\widehat{F_{1}}\sqsubseteq \widehat{F_{2}}$.\ On the other hand,\ suppose\ $\widehat{F_{1}}\sqsubseteq \widehat{F_{2}}$.\ We need to prove that\ $\widehat{F_{1}}\leq_P \widehat{F_{2}}$,\ that is to prove that\ $\{\widehat{F}\in PX: \widehat{F}\leq_{P} \widehat{F_{2}}\}$\ is a closed set in\ $PX$\ respect to\ $O_{\Rightarrow_P}(PX)$,\ by the proof of\ (1),\ $\{\widehat{F}\in PX: \widehat{F}\leq_{P}\widehat{F_{2}}\}\subseteq \overline{\{\widehat{F_{2}}\}}$,\ thus,\ when\ $\{\widehat{F}\in PX: \widehat{F}\leq_{P}\widehat{F_{2}}\}$\ is closed in\ $PX$\ respect to\ $O_{\Rightarrow_P}(PX)$,\ we have\ $\{\widehat{F}\in PX: \widehat{F}\leq_{P}\widehat{F_{2}}\}=\overline{\{\widehat{F_{2}}\}}$,\ it follows that\ $\widehat{F_{1}}\leq_P \widehat{F_{2}}$.\ Now,\  we prove that\ $\{\widehat{F}\in PX: \widehat{F}\leq_{P}\widehat{F_{2}}\}$\ is closed respect to\ $O_{\Rightarrow_P}(PX)$,\ equivalently,\ $\mathcal{U}=PX\setminus\{\widehat{F}\in PX: \widehat{F}\leq_{P}\widehat{F_{2}}\}$\ is a\ $\Rightarrow_P$\ convergence open set.\ Suppose\ $\mathcal{D}=\{\widehat{F_{i}}\}_{i\in I}\subseteq PX$\ is a\ $\leq_P$\ -\ directed set and\ $\mathcal{D}\Rightarrow_P \widehat{G}\in\mathcal{U}$.\ Let\ $G=\{a_1,a_2,\ldots,\ a_k\}$.\ By the definition of\ $\Rightarrow_P$\ convergence,\ there exists finite directed sets\ $D_{1},\dots,D_{n}\subseteq X$\ such that\ $D_1,\dots,D_n\subseteq\bigcup_{i\in I} \{ \da F_i\}$.\ For each\ $a_i\in G$,\ there is one of the above directed set\ $D_i$\ such that\ $D_i\rightarrow a_i$,\ $i=1,2,\ldots,k$.\ By contradiction,\ suppose\ ${\mathcal D}\cap {\mathcal U}=\emptyset$,\ then,\ by the definition of order\ $\leq_{P}$,\ $\forall \widehat{F_i}\in\mathcal{D},\ \da F_i\subseteq\da F_2$,\  thus\ $\bigcup {\da F_i}\subseteq \da F_2$,\ it follows that\ $D_i\subseteq \da F_2,\ i=1,2,\ldots,\ k$.\ But\ $\da F_2$\ is a closed set of\ $X$,\ then,\ all limit points of\ $D_i$\ belongs to\ $\da F_2$,\ it follows that\ $G=\{a_1,a_2,\ \ldots,\ a_k\}\subseteq \da F_2$.\ However,\ according to the hypothesis,\ $\widehat{G}\in\mathcal{U}$,\ that is $\widehat{G}\nleq_P\widehat{F_2}$,\ then,\ it follows that\ $\ua F_2\nsubseteq \ua G$,\  that is\ $\exists a\in F_2$\ such that\ $G\subseteq X\setminus\da a$.\ By the definition of\ $\Rightarrow_P$\ convergence,\ each\ $D_i$\ converges to the points in\ $G$,\ and\ $X\setminus\da a$ is an open set of $X$,\ then\ $D_i\cap(X\setminus\da a)\neq\emptyset$.\ For each\ $i=1,\dots,n$,\ pick\ $d_i\in D_i\cap(X\setminus\da a)$.\ We have $(d_1,\dots,d_n)\in\prod_{i=1}^n D_i$\ and\ $a\notin\ua (d_1,\dots,d_n)$.\ By the hypothesis,\ $\mathcal{D}\cap\mathcal{U}=\emptyset$,\ it follows that\ $\forall\widehat{F_i}\in \mathcal{D},\ \ua F_2\subseteq\ua F_i $.\ Then\ $\forall i\in I,\ a\in\ua F_i,\ua F_i\nsubseteq \ua(d_1,\dots,d_n)$.\ This contradicts with\ (4)\ in the definition of\ $\Rightarrow_P$\ convergence.\ Therefore,\ $\mathcal{U}$\ is a\ $\Rightarrow_P$\ convergence open set in\ $PX$.

\vskip 3mm

(3) For an arbitrary topological space\ $X$,\ $O(X)\subseteq d(X)$\ holds,\ then\ $O_{\Rightarrow_P}(PX)\subseteq d(PX)$.\ On the other hand,\ according to the definition of\ $\Rightarrow_P$\ convergence topology,\ if directed set\ $\mathcal D\in PX$\ with\ $\mathcal{D}\Rightarrow_{P}\widehat{F}$,\ then\ $\mathcal{D}$\ converges to\ $\widehat{F}$\ respect to\ $O_{\Rightarrow_P}(PX)$.\ Thus,\ by the definition of directed open set,\ $\mathcal{D}\Rightarrow_{P}\widehat{F} \in\mathcal{U}\in d(PX)$\ implies\ $\mathcal{U}\cap \mathcal{D}\neq \emptyset$.\ Then,\ $\mathcal{U}\in O_{\Rightarrow_P}(PX)$,\ it follows that\ $O_{\Rightarrow_P}(PX)= d(PX)$,\ that is\ $(PX,\  O_{\Rightarrow_P}(PX))$ is a directed space.\ $\Box$

\begin{proposition}\label{prop PXdc}
 Suppose\ $X,\ Y$\ are two directed spaces.\ Then function\ $f:(PX,\ O_{\Rightarrow_P}(PX))\rightarrow Y$\ is continuous if and only if for each directed set\ ${\mathcal D}\subseteq PX$\ and\ $\widehat{F}\in PX$,\ ${\mathcal D}\Rightarrow_P\widehat{F}$\ implies\ $f(\mathcal{D})\rightarrow f(\widehat{F})$.
\end{proposition}

\noindent{\bf Proof} Since\ $\Rightarrow_P$\ convergence will lead to\ $O_{\Rightarrow_P}(PX)$\ topological convergence,\ the necessity is obvious.\ We are going to prove the sufficiency.\ Firstly,\ we check that\ $f$\ is monotone.\ If\ $\widehat{F_1},\ \widehat{F_2}\in PX$\ and\ $\widehat{F_1}\leq_P\widehat{F_2}$,\ then\ $\{\widehat{F_2}\}\Rightarrow_P\widehat{F_1}$.\ By the hypothesis,\ $\{f(\widehat{F_2})\}\rightarrow f(\widehat{F_1})$,\ it follows that\ $f(\widehat{F_2})\sqsubseteq f(\widehat{F_1})$.\ Suppose\ $U$\ is an open set of\ $Y$\ and there is a directed set\ ${\mathcal D}\Rightarrow_P\widehat{F}\in f^{-1}(U)$,\ then\ $f(\mathcal{D})$\ is a directed set of\ $Y$\ and\ $f(\mathcal{D})\rightarrow f(\widehat{F})\in U$,\ thus,\ there exists some\ $ \widehat{F'}\in {\mathcal D}$\ such that\ $f(\widehat{F'})\in U$,\ that is,\ $\widehat{F'}\in {\mathcal D}\cap f^{-1}(U)$.\ According to the definition of\ $\Rightarrow_P$\ convergence open set,\ $f^{-1}(U)\in O_{\Rightarrow_P}(PX)$,\ that is\ $f$\ is continuous.\ $\Box$

\vskip 3mm

Define a binary operation\ $\oplus:PX\times PX\rightarrow PX$\ on\ $PX$\ :$\forall\ \widehat{F_1},\ \widehat{F_2}\in PX,\  \widehat{F_1}\oplus\widehat{F_2}=\widehat{F_1\cup F_2}$.\ For arbitrary\ $\widehat{F_1}=\widehat{F_2},\ \widehat{G_1}=\widehat{G_2}$\ with\ $F_1\neq F_2,\ G_1\neq G_2$.\ $F_1\cup G_1\subseteq\da(F_2\cup G_2)$\ and\ $F_1\cup G_1\subseteq\ua(F_2\cup G_2)$\ imply\ $\widehat{F_1\cup G_1}\leq_P\widehat{F_2\cup G_2}$.\ Similarly,\ we have the opposite inequality,\ thus,\ $\oplus$\ is well-defined.

\vskip 3mm

\begin{theorem} Let\ $X$\ be a directed space.\ Then\ $(PX,\ \oplus)$\ is a directed semilattice.
\end{theorem}

\noindent{\bf Proof} By Proposition\ \ref{prop PX},\ $(PX,\ O_{\Rightarrow_P}(PX))$\ is a directed space.\ We will prove that\ $\oplus $\ is a directed semilattice operation.\ It is evidently to check that\ $\oplus$\ satisfy \ (a),\   (b),\   (c)\ in definition\ \ref{def dsl},\ we now prove the continuity of\ $\oplus$.
$\oplus$\ is obviously monotone.\ By Theorem\ \ref{theorem opc}(2)\ and proposition\ \ref{prop PXdc},\ we only need to prove that,\ for each directed set\ ${\mathcal D}=\widehat{F_i}\subseteq PX$\ and\ $\widehat{F},\ \widehat{G}\in PX$,\ ${\mathcal D}\Rightarrow_P\widehat{F}$\ will imply\ $\widehat{G}\oplus {\mathcal D}\Rightarrow_P(\widehat{G}\oplus \widehat{F})=\widehat{(G\cup F)}$.\ Here,\ $G\cup {\mathcal D}=\{\widehat{G}\oplus \widehat{F_i}:\widehat{F_i}\in \mathcal{D}\}$\ is also a directed set.\ By the definition of\ $\Rightarrow_P$\ convergence,\ there exists finite directed sets\ $D_1,\dots,\ D_n\subseteq X$\ such that\ $\mathcal{D}\Rightarrow_P\widehat{F}$.\ Let\ $G=\{a_1,\dots,a_k\}$,\ take\ $D_{n+1}=\{a_1\},\dots,D_{n+k}=\{a_k\}$.We can evidently verify that\ $D_1,\dots,D_{n+k}$\ is sufficient to satisfy the definition of\ $\widehat{G}\oplus\mathcal{D}\Rightarrow_P\widehat{G\cup F}$.\ Therefore,\ $(PX,\ \oplus)$\ is a directed semilattice.\ $\Box$

\vskip 3mm

The following theorem is the main result of this section.

\begin{theorem}\label{theorem PX}
	Suppose\ $X$\ is a directed space,\ then\ $(PX,\ O_{\Rightarrow_P}(PX))$\ is the directed convex powerspace over\ $X$,\ that is,\ endowed with topology\ $O_{\Rightarrow_P}(PX)$,\ $(PX,\ \oplus)\cong P_P(X)$.	
\end{theorem}

\noindent{\bf Proof} Define function\ $i:X\rightarrow PX$\ as follows:\ $\forall x\in X$,\ $i(x)=\widehat{x}$.\ We prove the continuity of\ $i$.\ $i$\ is evidently monotone.\ Suppose we have directed set\ $D\subseteq X$\ and\ $x\in X$\ satisfy\ $D\rightarrow x$.\ Let\ $\mathcal{D}=i(D)=\{\widehat{d}:d\in D\}$,\ then\ $\mathcal{D}$\ is a directed set in\ $PX$\ and\ ${\mathcal D}\Rightarrow_{P}\widehat x$.\ Since\ $D$\ satisfy the definition such that\ $\mathcal D\Rightarrow_P\widehat{x}$.\ By Proposition\ \ref{prop dc},\ $i$\ is continuous.

Let\ $(Y,\   +)$\ be a directed semilattice,\ $f:X\rightarrow Y$\ is a continuous function.\ Define\ $\bar{f}:PX\rightarrow Y$\ as follows:\ $\forall \widehat{F}\in PX$ (let\ $F=\{a_1,a_2,\ldots,a_n\}$),\ $$\bar{f}(\widehat{F})=f(a_1)+ f(a_2)+\dots+ f(a_n)=\sum\limits_{a\in F}f(a)$$
Particularly,\ we write\ $\bar{f}(\widehat{F})=\Sigma f(F)$.\\
$\bar{f}$\ is well-defined,\ since if we have\ $\widehat{F_1},\ \widehat{F_2}\in PX$\ and\ $\widehat{F_1}=\widehat{F_2}$\ but with\ $F_1\neq F_2$.\ Let\ $F_1=\{a_1,\dots,a_k\},\ F_2=\{b_1,\dots,b_n\}$.\ For\ $\da F_1\subseteq\da F_2$,\ we have\ $b_1,\dots,b_k\in F_2$\ such that\ $a_1\leq b_1,\dots,a_k\leq b_k$.\ If\ $F_2\setminus \{b_1,\dots,b_k\}=\emptyset$,\ then\ $\Sigma f(F_1)=\bar{f}(\widehat{F_1})\leq\Sigma f(F_2)=\bar{f}(\widehat{F_1})$.\ If not,\ suppose\ $F_2\setminus\{b_1,\dots,b_k\}=\{b_{i1},\dots,b_{is}\}$,\ but\ $\ua F_2\subseteq\ua F_1$,\ we have\ $a_{i1},\dots,a_{is}\in F_1$\ such that\ $a_{i1}\leq b_{i1},\ \dots,\ a_{is}\leq b_{is}$,\ now,\ we have\ $a_1\leq b_1,\dots,a_n\leq b_n,\ a_{i1}\leq b_{i1},\dots,a_{is}\leq b_{is}$.\ Notice that\ $a_{i1},\dots,a_{is}$\ repeat\ $\{a_1,\dots,a_n\}$,\ it follows that\ $\bar{f}(\widehat{F_1})= f(a_1)+\dots +f(a_n)+f(a_{i1})+\dots +f(a_{is})\leq f(b_1)+\dots +f(b_n)+f(b_{i1})+\dots +f(b_{is})=\bar{f}(\widehat{F_2})$.\ Similarly,\ we have\ $\bar{f}(\widehat{F_2})\leq \bar{f}(\widehat{F_1})$.\ Thus,\ $\bar{f}$ is well-defined.

(1) $f=\bar{f}\circ i$.

For arbitrary\ $x\in X$,\ $(\bar{f}\circ i)(x)=\bar{f}(i(x))=\bar{f}(\widehat{x})=f(x)$.

(2) $\bar{f}$\ is a directed semilattice homomorphism,\ that is,\ $\bar{f}$\ is continuous and for arbitrary\ $\widehat{F_1},\   \widehat{F_2}\in PX$,\ $\bar{f}(\widehat{F_1}\oplus \widehat{F_2})=\bar{f}(\widehat{F_1})+ \bar{f}(\widehat{F_2})$.

Firstly,\ we prove that\ $\bar{f}$\ preserves directed semilattice operation.\ Suppose\ $\widehat{F_1},\ \widehat{F_2}\in PX$.\ Then\ $\bar{f}(\widehat{F_1}\oplus\widehat{F_2})=\bar{f}(\widehat{F_1\cup F_2})=\Sigma f(F_1\cup F_2)=(\Sigma f(F_1)+(\Sigma f(F_2))=\bar{f}(\widehat{F_1})+\bar{f}(\widehat{F_2})$.\ Next,\ we prove the continuity of\ $\bar{f}$.\ We check that\ $\bar{f}$\ is monotone.\ Suppose\ $\widehat{F_1},\ \widehat{F_2}\in PX$.\ Let\ $F_1=\{a_1,\dots,a_k\},\ F_2=\{b_1,\dots,\ b_n\}$,\ and\ $\widehat{F_1}\leq_P\widehat{F_2}$,\ that is\ $\da F_1\subseteq\da F_2\ \& \ua F_2\subseteq\ua F_1$.\ For\ $\da F_1\subseteq\da F_2$,\ we have\ $b_1,\dots,b_k\in F_2$\ such that\ $a_1\leq b_1,\dots,a_k\leq b_k$.\ If\ $F_2\setminus \{b_1,\dots,b_k\}=\emptyset$,\ then\ $\Sigma f(F_1)=\bar{f}(\widehat{F_1})\leq\Sigma f(F_2)=\bar{f}(\widehat{F_1})$.\ If not,\ suppose\ $F_2\setminus\{b_1,\dots,b_k\}=\{b_{i1},\dots,b_{is}\}$,\ but\ $\ua F_2\subseteq\ua F_1$,\ we have\ $a_{i1},\dots,a_{is}\in F_1$\ such that\ $a_{i1}\leq b_{i1},\ \dots,a_{is}\leq b_{is}$,\ now,\ we have\ $a_1\leq b_1,\dots,a_n\leq b_n,\ a_{i1}\leq b_{i1},\dots,\ a_{is}\leq b_{is}$.\ Notice that\ $a_{i1},\ \dots,\ a_{is}$\ repeat\ $\{a_1,\ \dots,\ a_n\}$,\ it follows that\ $\bar{f}(\widehat{F_1})= f(a_1)+\dots +f(a_n)+f(a_{i1})+\dots +f(a_{is})\leq f(b_1)+\dots +f(b_n)+f(b_{i1})+\dots +f(b_{is})=\bar{f}(\widehat{F_2})$.\ In conclusion,\ $\widehat{F_1}\leq_P\widehat{F_2}$\ will imply\ $\bar{f}(\widehat{F_1})\leq \bar{f}(\widehat{F_2})$,\ that is,\ $\bar{f}$\ is monotone.

Finally,\ we verify the continuity of\ $\bar{f}$.\ Suppose\ ${\mathcal D}\subseteq PX$\ is a directed set and\ ${\mathcal D}\Rightarrow_P\widehat{F}\in PX$,\ there exists finite directed sets\ $D_{1},\dots,\  D_{n}\subseteq X$\ satisfy definition of\ $\Rightarrow_P$\ convergence.\ Let\ $F=\{b_1,b_2,\ldots,b_k\}$.\ By\ (2)\ in the definition,\ for each\ $D_i,\ 1\leq i\leq n$,\ there exists some\ $b_i\in F$\ such that\ $D_i\rightarrow b_i$.\ If\ $F\setminus \{b_1,\dots,b_n\}\neq\emptyset$,\ which is denoted by\ $G=\{a_1,a_2\ldots,a_s\}$.\ By\ (3)\ of definition,\ for each\ $a_j\in G$,\ there exists\ $1\leq i_j\leq n$\ such that\ $D_{i_j}\rightarrow a_j$.\ By the continuity of\ $f$,\ we have\ $f(D_{i})\rightarrow f(b_{i}),\ i=1,\dots,n$,\ and\ $f(D_{i_j})\rightarrow f(a_j)$,\ $j=1,2,\ldots,s$.\ Since the directed semilattice operation\ $+$\ on\ $Y$\ is continuous,\ the following convergence holds
$$f(D_{1})+\cdots + f(D_n)+ f(D_{i_1})+\cdots+ f(D_{i_s})\rightarrow f(b_{1})+ \dots + f(b_n)+ f(a_{i_1})+\cdots+ f(a_{i_s}).\ \\ (\ast)$$
Here,\ $f(D_{1})+ \cdots + f(D_n)+  f(D_{i_1})+\cdots+ f(D_{i_s})=\{f(d_{1}) + \cdots + f(d_n)+ f(d_{i_1})+\cdots+ f(d_{i_s}) :(d_1,\ldots,d_k,d_{i_1},\ldots,\ d_{i_s})\in (\prod\limits_{i=1}^{n}D_i\})\times(\prod\limits_{j=1}^{s}D_{i_j})\}$.
Let\ $U$\ be an arbitrary open neighborhood of\ $\Sigma f(F)$,\ by\ $(\ast)$,\ there exists some\ $(d_1,\ldots,d_n,d_{i_1},\ldots,d_{i_s})\in (\prod\limits_{i=1}^{n}D_i\})\times(\prod\limits_{j=1}^{s}D_{i_j})$\ such that\ $f(d_{1})+
 \cdots + f(d_{n})+ f(d_{i_1})+\cdots+ f(d_{i_s})\in U$.\ Since each\ $D_{i_j}$\ repeats\ $D_i$,\ there exists\ $(d^{'}_1,d^{'}_2,\ldots,d^{'}_n)\in\prod\limits_{i=1}^{n}D_i\}$ such that\ $f(d^{'}_1)+\cdots+ f(d^{'}_n)\sqsupseteq f(d_{1})+
 \cdots + f(d_{n})+ f(d_{i_1})+\cdots+ f(d_{i_s})$.\ By\ (4)\ in the definition,\ there exists some\ $\widehat{F_{1}}\in\mathcal{D}$\ such that\ $\uparrow F_{1}\subseteq
 \uparrow (d^{'}_{1},\ \dots,\ d^{'}_{n})$.\ On the other hand,\ by\ (1)\ in definition of\ $\Rightarrow_{P}$\ convergence,\ for each\ $d_i^{'}$,\ there exists some\ $\widehat{F_i}\in\mathcal{D}$\ such that\ $d^{'}_i\in\da F_i$.\ Since\ $\mathcal D$ is directed,\ we have\ $\widehat{F_{2}}\in\mathcal{D}$\ with\ $\widehat{F_i}\leq_P\widehat{F_{2}},\ i=1,\ \dots,\ n$.\ Again,\ by the directness of\ $\mathcal D$,\ we obtain\ $\widehat{F'}\in \mathcal D$\ such that $\widehat{F_1},\ \widehat{F_2}\leq_P\widehat{F'}$.\ Now,\  we have\ $\da(d^{'}_1,\dots,d^{'}_n)\subseteq \da F_2\subseteq \da F^{'},\ \ua F^{'}\subseteq\ua F_1\subseteq\ua (d^{'}_1,\dots,d^{'}_n)$,\ that is $\widehat{\{d^{'}_1,\dots,d^{'}_n\}}\leq_P\widehat{F^{'}}$,\ since\ $\bar{f}$\ is monotone and\ $U$\ is an upper set in\ $Y\ $,\ we got\ $\bar{f}(\widehat{F^{'}})\in U$.\ It follows that\ $\bar{f}({\mathcal D})=\{\Sigma f(F_i):\widehat{F_i}\in\mathcal{D}\}\rightarrow\Sigma f(F)$.\ By Proposition\ \ref{prop PXdc},\ $\bar{f}$\ is continuous.

(3) \ $\bar{f}$\ is unique.

Suppose we have a directed semilattice homomorphism\ $g:(PX,\ \oplus)\rightarrow (Y,\ +)$\ such that\ $f=g\circ i$.\ Then\ $g(\widehat{x})=f(x)=\bar{f}(\widehat{x})$.\ For each\ $ \widehat{F}\in PX$\ (Let\ $ F=\{a_{1},\   \dots ,\   a_{n}\}$\ ),\ \begin{eqnarray*}
 g(\widehat{F})&=& g(\widehat{a_1}\oplus\widehat{a_2}\oplus\cdots\oplus\widehat{a_n})\\ &=&g(\widehat{a_{1}})\oplus g(\widehat{a_2})\oplus \cdots \oplus g(\widehat{a_{n}})\\ &=&\bar{f}(\widehat {a_{1}})\oplus \bar{f}(\widehat{a_2})\oplus\cdots\oplus\bar{f}(\widehat{a_{n}})\\ &=&\bar{f}(\widehat{a_1}\oplus\widehat{a_2}\oplus\cdots\oplus\widehat{a_n})\\ &=&\bar{f}(\widehat F).
 \end{eqnarray*}
 That is\ $\bar{f}$\ is unique.

 In conclusion,\ by Definition\ \ref{def dcps},\ endowed with the topology\ $O_{\Rightarrow_P}(PX)$,\ the directed semilattice\ $(PX,\ \oplus)$\ is a directed convex powerspace over\ $X$,\ that is,\ $P_P(X)\cong (PX,\ \oplus)$.\ $\Box$

 \vskip 3mm

 The directed convex powerspace is unique in the sense of order isomorphism and topological
homomorphism,\ so we can denote the directed convex powerspace of each\ $X$\ by\ $P_P(X)=(PX,\ \oplus)$.

 Suppose\ $X,\ Y$\ are two directed spaces,\ $f:X\rightarrow Y$\ is a continuous function.\ Define map\ $P_P(f):P_P(X)\rightarrow P_P(Y)$\ as follows:\ $\forall \widehat{F}\in PX$,\ $$P_P(f)(\da F)=\widehat{f(F)}.$$
It is evident that,\ $P_P(f)$\ is well-defined and order preserving.\ According to the proof of the theorem above,\  it is easy to check that\ $P_P(f)$\ is a directed semilattice homomorphism
between these two directed convex powerspaces.\ If\ $id_X$\ is the identity function and\ $g:Y\rightarrow Z$\ is an arbitrary continuous function from\ $Y$\ to a directed space\ $Z$,\ then,\ $P_P(id_X)=id_{P_P(X)},\ P_P(g\circ f)=P_P(g)\circ P_P(f)$.\ Thus,\ $P_P:{\bf Dtop}\rightarrow {\bf Dsl}$\ is a functor from\ ${\bf Dtop}$\ to\ ${\bf Dsl}$.\ Let\ $U:{\bf Dsl}\rightarrow {\bf Dtop}$\ be the forgetful functor.\ By Theorem\ \ref{theorem PX},\ we have the following results.

 \begin{corollary}
$P_P$\ is a left adjoint of the forgetful functor\ $U$,\ that is,\ ${\bf Dsl}$\ is a reflective subcategory of \ ${\bf Dtop}$.
 \end{corollary}

%%%%%%%%%%%%%%%%%%%%%%%%%%%%%%%%%%%%%%%%%%%%%%%%%%%%%%%%%%%%%%

\section{Relations Between Convex Powerspaces}

In this section,\ we will discuss the relations between the convex powerdomain of
dcpo and its directed lower powerspace.

Let\ $L$\ be a domain equipped with the\ Scott\ topology.\ A nonempty subset\ $A$\ is a\ $lens$\ if\ $A$\ can be written as the intersection of a closed set and a compact saturated set.
 A lens\ $A=C\cap K$\ has a canonical representation of the form\ $\bar A\cap\ua A$\ (since\ $A\subseteq \bar A\cap\ua A\subseteq C\cap K=A$)\ ;note that \ $\ua A$\ is compact since\ $A$\ is.A pair\ $(C,\ K)\in C(L)\times Q(L)$\ is called a lens factorization if\ $C=\bar A$\ and\ $K=\ua A$\ for some lens\ $A$.We denote by\ $Lens L$ the set of all lens.

\begin{theorem}{\rm\cite{GHK}} Let\ $L$\ be a domain equipped with the\ Scott\ topology.\ If\ $L$\ is coherent,\ resp.\ countably based,\ then the convex powerdomain\ $P(L)$\ (which is also called the Plotkin powerdomain)may be identified with Lens\ $L$\ equipped with the Egli-Milner order,\ resp.\ topological Egli-Milner order.

\end{theorem}

\begin{proposition}
Let\ $X$\ be a continuous domain endowed with\ Scott\ topology,\ then\ $\sigma(P(X))|_{PX}=O_{\Rightarrow_P}(PX)$.

\end{proposition}
\noindent{\bf Proof}  By proposition\ \ref{prop rstLX}\ and proposition\ \ref{prop rstUX}(1),\ it is not hard to check that\ $\sigma(P(X))|_{PX}\subseteq O_{\Rightarrow_P}(PX)$.

On the other hand,\ for each\ $\mathcal U\in O_{\Rightarrow_P}(PX)$,\ let$$\ua_{P(X)}\mathcal U =\{(C,\ D)\in C(X)\times Q(X):\exists \widehat{F}\in PX \ such\ that\ \da F\subseteq C \&D\subseteq\ua F\}.$$ Then,\ it is obviously that\ $\mathcal U=\ua_{P(X)}\mathcal U \cap PX$,\ it is sufficient to prove that\ $\ua_{P(X)}\mathcal U\in\sigma(P(X))$.\ Let\ $\mathcal D=\{(C_i,\ D_i)\}_{i\in I}\subseteq P(X)$\ be a directed set with\ $\vee\mathcal D=(\overline{\cup_{i\in I}C_i},\ \cap_{i\in I}D_i)\in \ua_{P(X)}\mathcal U$.\ By the definition of\ $\ua_{P(X)}\mathcal U$\ we have some\ $\widehat{F}\in\mathcal U$\ such that\ $\da F\subseteq C\ \&\ D\subseteq\ua F$.\ Set\ $F=\{a_1,\ \dots,\ a_n\}$.\ Let\ $D_i=\dda a_i,\ 1\leq i\leq n$,\ $\mathcal F=\{\widehat{\{d_1,\ \dots ,\ d_n\}}:(d_1,\ \dots ,\ d_n)\in\prod\limits_{i=1}^n D_i\}$.\ Since\ $X$\ is a continuous domain and by the proof of proposition\ \ref{prop rstUX}(1),\ we can directly check that\ $\mathcal F\Rightarrow_P\widehat{F}$,\ it follows that there exists some\ $\widehat{\{d_1,\ \dots,\ d_n\}}\in\mathcal F\cap\mathcal U$.\ Since\ $F\subseteq\overline{\cup_{i\in I}C_i}$\ and\ $X$\ is continuous,\ we can infer that each\ $D_i\subseteq\da\cup_{i\in I}C_i=\cup_{i\in I}C_i$.\ For\ $\mathcal D$\ is directed,\ we have some\ $(C^{'},\ D^{'})\in\mathcal D$\ such that\ $\da\{d_1,\ \dots ,\ d_n\}\subseteq C^{'}$.\ Besides,\ $\cap_{i\in I}D_i\subseteq\ua F\subseteq(\cup_{i=1}^n\ua d_i)^{\circ}$,\ since\ $X$\ is a continuous domain endowed with the Scott topology,\  hence well-filtered,\  we have some\ $(C^{''},\ D^{''})\in\mathcal D$\ such that\ $D^{''}\subseteq(\cup_{i=1}^n\ua d_i)^{\circ}\subseteq\ua\{d_1,\ \dots ,\ d_n\}$.\ Again,\ by the directness of\ $\mathcal D$,\ choose\ $(C,\ D)\in\mathcal D$\ with\ $C^{'},\ C^{''}\subseteq C\ \&D\subseteq D^{'},\ D^{''}$.\ Then,\ $\da\{d_1,\ \dots,\ d_n\}\subseteq C^{'}\subseteq C\ \&D\subseteq D^{''}\subseteq\ua\{d_1,\ \dots,\ d_n\}$,\ that is\ $(C,\ D)\in\ua_{P(X)}\mathcal U $\ and\ $\ua_{P(X)}\mathcal U $\ is a\ Scott\ open set.$\Box$

\begin{example}
Let\ $X$ be the same as in example\ \ref{example UXneq},\ by example \ref{example UXneq}\ we can directly infer that,\ $P_P(P_O(X))\neq P(P_O(X))$,\ that is,\ in general,\ the directed convex powerspace over a\ dcpo\ endowed with the\ Scott\ topology is not agree with its convex powerdomain.

\end{example}

\section{The Commutativity of The Directed Upper and Lower Functors}

In domain theory,\ the Hoare and Smyth powerdomain constructors were firtsly proved to be commuted under composition in\ \cite{FM90}\ in 1990.\ In 1991,\ Heckmann gave an algebraic methods and does not rely on any explicit representations of the powerdomains\ \cite{HECK91}.\ In this section,\ we will discuss the commutativity of the directed upper and lower functors.

In this paper,\ for convenience,\ a directed inflationary semilattice is abbriated as\ $\vee-$\ directed space,\ a directed deflationary semilattice is abbriated as\ $\wedge-$\ directed space.\ A morphism between two directed deflationary semilattices is called a\ $\vee-$\ morphism,\ a morphism between two directed inflationary semilattices is called a\ $\wedge-$\ morphism.

\begin{definition}\
	\begin{enumerate}[(1)]
		\item 	A\ $\vee-\wedge-$\ directed space\ $X$\ is a directed space that is both a\ $\vee-$\ directed space and\ $\wedge-$\ directed space that the distributive law\ $a\wedge(b\vee c)=(a\wedge b)\vee(a\wedge c),\ \forall a,b,c\in X$.
		\item A\ $\vee-\wedge-$\ morphism between two\ $\vee-\wedge-$\ directed spaces that is both a\ $\vee-$\ morphism and a\ $\wedge-$\ morphism.
	\end{enumerate}

\end{definition}
 Denote the category of all\ $\vee-\wedge-$\ directed spaces together with\ $\vee-\wedge-$\ morphisms by\ ${\bf Ddl}$.\ For all directed spaces\ $X$,\ we then show that both\ $P_U(P_L(X))$\ and\ $P_L(P_U(X))$\ are free\ $\vee-\wedge-$\ directed spaces over\ $X$,\ whence they are isomorphic.

We will discuss some properties of functor\ $P_L$,\ and functor\ $P_U$\ is analogously.\ Theorm\ \ref{theorem LX}\ indicates that\ for each directed space\ $X$,\ $P_L(X)$\ is a\ $\vee-$\ directed space.\ The following proposition tell us that if\ $X$\ is a\ $\wedge-$\ directed space,\ then\ $P_L(X)$\ is actually a\ $\vee-\wedge-$\ directed space.

\begin{theorem}\label{comu 1}\
	\begin{enumerate}[(1)]
		\item If\ $X$\ is a\ $\wedge-$\ directed space,\ then\ $P_L(X)$\ is a\ $\vee-\wedge-$\ directed space,\ and\ $\eta_{\vee}:X\rightarrow P_L(X)$\ is a\ $\wedge-$\ morphism,\ here\ $\eta_{\wedge}(x)=\da x,\forall x\in X$.
		\item  If\ $X$\ is a\ $\vee-$\ directed space,\ then\ $P_U(X)$\ is a\ $\vee-\wedge-$\ directed space,\ and\ $\eta_{\wedge}:X\rightarrow P_U(X)$\ is a\ $\vee-$\ morphism,\ here\ $\eta_{\wedge}(x)=\uparrow x,\forall x\in X$.
	\end{enumerate}
\end{theorem}

{\bf Proof}\ \ We only prove\ (1).\ Let\ $X$\ be a\ $\wedge-$\ directed space,\ then\ $P_L(X)$\ is a\ $\vee-$\ directed space,\ we have to construct the operation\ $\wedge$\ in\ $P_L(X)$.\ To this end,\ for each\ $\da F_1,\da F_2\in P_L(X)$,\ let\ $f_a(x)=\eta_{\vee}(a\wedge x)$\ for each\ $a\in X$.\ By Theorem\ \ref{theorem LX},\ there exist a unique morphism\ $\bar{f_a}:P_L(X)\rightarrow P_L(X)$\ such that\ $\bar{f_a}(\da F)=\da f_a(F)$\ for each\ $\da F\in P_L(X)$.\ Let\ $g(a)=\bar{f_a}(\da F_2)=\da\{a\wedge b:b\in F_2\}$,\ again by Theorem\ \ref{theorem LX},\ there exist a unique morphism\ $\bar{g}:P_L(X)\rightarrow P_L(X)$\ such that\ $\bar{g}(\da F)=\da\{a\wedge b:a\in F,b\in F_2\}$.\ Thus,\ we have a continuous operation in\ $P_L(X)$\ defined as follows
$$\da F_{1}\sqcap\da F_2=\{\eta_{\wedge}(a\wedge b):a\in F_{1},b\in F_{2}\}=\da\{a\wedge b:a\in F_{1},b\in F_{2}\}$$

It is straightly to check that\ $\sqcap$\ satisfy all the conditions in Definition\ \ref{def disl}:

\begin{enumerate}[(1)]
	\item $\da F\sqcap\da F=\da\{a\wedge a:a\in F\}=\da F,\forall\da F\in P_L(X)$.\\
	\item $(\da F_1\sqcap\da F_2)\sqcap\da F_3=\da\{a\wedge b:a\in F_1,b\in F_2\}\sqcap\da F_3=\da\{a\wedge b\wedge c:a\in F_1,b\in F_2,c\in F_3\}=\da F_1\sqcap\da\{b\wedge c:b\in F_2,c\in F_3\}=\da F_1\sqcap(\da F_2\sqcap\da F_3),\forall\da F_1,\da F_2,\da F_3\in P_L(X)$.\\
	\item $\da F_1\sqcap\da F_2=\da\{a\wedge b:a\in F_1,b\in F_2\}=\da F_2\sqcap\da F_1,\forall\da F_1,\da F_2\in P_L(X)$.\\
	\item $\da F_1\sqcap\da F_2=\da\{a\wedge b:a\in F_1,b\in F_2\}\leq\da F_1,\forall\da F_1,\da F_2\in P_L(X)$.
\end{enumerate}

Next,\ we shall check that the distributivity law for\ $\cup$\ and\ $\sqcap$.\ Suppose we have\ $F_1=\{a_1,\dots,a_n\},\ F_2=\{b_1,\dots,b_m\},\ F_3=\{c_1,\dots,c_k\}$.

\begin{enumerate}[(1)]
	\item $\da F_1\sqcap(\da F_2\cup\da F_3)=(\da F_1\sqcap\da F_2)\cup(\da F_1\sqcap\da F_3)$.\\
	$\da F_1\sqcap(\da F_2\cup\da F_3)=\da F_1\sqcap\da\{b_1,\dots,b_m,c_1,\dots,b_k\}=\da\{a\wedge b,a\wedge c:a\in F_1,b\in F_2,c\in F_3\}$.\\
	$(\da F_1\sqcap\da F_2)\cup(\da F_1\sqcap\da F_3)=\da\{a\wedge b:a\in F_1,b\in F_2\}\cup\da\{a\wedge c:a\in F_1,c\in F_3\}=\da\{a\wedge b,a\wedge c:a\in F_1,b\in F_2,c\in F_3\}$.\ Thus they are equal.\\
	
	\item $\da F_1\cup(\da F_2\sqcap\da F_3)=(\da F_1\cup\da F_2)\sqcap(\da F_1\cup\da F_3)$.\\
	$\da F_1\cup(\da F_2\sqcap\da F_3)=\da F_1\cup\da\{b\wedge c:b\in F_2,c\in F_3\}=\da\{a,b\wedge c:a\in F_1,b\in F_2,c\in F_3\}$.
	$(\da F_1\cup\da F_2)\sqcap(\da F_1\cup\da F_3)=\da\{a\wedge a,a\wedge c,a\wedge b,b\wedge c:a\in F_1,b\in F_2,c\in F_3\}=\da\{a,b\wedge c:a\in F_1,b\in F_2,c\in F_3\}$.\ Thus they are equal.
	
\end{enumerate}
Hence,\ $(P_L(X),\vee,\sqcap)$\ is a\ $\vee-\wedge-$\ directed space.\ Last,\ we shall prove that\ $\eta_\vee$\ is a\ $\wedge-$\ morphism.\ For each\ $x,\ y\in X,\eta_\vee(x\wedge y)=\da(x\wedge y)=(\da x)\sqcap(\da y)=\eta_\vee(x)\sqcap\eta_\vee(y)$.\ $\Box$

In the pervious part,\ we showed that the functor\ $P_L$ does not destroy the\ $\wedge$\ operator,\ now we prove that\ $P_L$\ has the analogous property for\ $\wedge-$\ morphism.

\begin{theorem}\label{comu 2}\
	\begin{enumerate}[(1)]
		\item Let\ $X$\ be a\ $\wedge-$\ directed space,\ $Y$\ a\ $\vee-\wedge-$\ directed space,\ and\ $f:X\rightarrow Y$\ a\ $\wedge-$\ morphism.\ Then its unique extension\ $P_L(f):P_L(X)\rightarrow Y$\ is a\ $\vee-\wedge-$\ morphism.
		\item Let\ $X$\ be a\ $\vee-$\ directed space,\ $Y$\ a\ $\vee-\wedge-$\ directed space,\ and\ $f:X\rightarrow Y$\ a\ $\vee-$\ morphism.\ Then its unique extension\ $P_U(f):P_U(X)\rightarrow Y$\ is a\ $\vee-\wedge-$\ morphism.
	\end{enumerate}
	
\end{theorem}
{\bf Proof}\ \ We only prove\ (1).\ Only need to check that\ $P_L(f)$\ is a\ $\wedge-$\ morphism.\ Suppose we have\ $F_1={a_1,\dots,a_n},\ F_2=\{b_1,\dots,b_m\}$.\ \begin{eqnarray*} P_L(f)(\da F_1\sqcap\da F_2)&=&P_L(f)(\da\{a\wedge b:a\in F_1,b\in F_2\})\\ &=&\bigvee\limits_{a\in F_1,b\in F_2}f(a\wedge b)\\ &=&\bigvee\limits_{a\in F_1,b\in F_2}f(a)\wedge f(b)
\end{eqnarray*}

On the other hand,\ \begin{eqnarray*} P_L(f)(\da F_1)\wedge P_L(f)(\da F_2)&=&(\vee f(F_1))\wedge (\vee f(F_2))\\ &=&(\vee f(F_1)\wedge f(b_1))\vee (\vee f(F_1)\wedge f(b_2))\vee \dots \vee (\vee f(F_1)\wedge f(b_m))\\ &=&(f(a_1)\wedge f(b_1))\vee \dots \vee (f(a_n)\wedge f(b_1))\vee \dots \vee (f(a_n)\wedge f(b_m))\\ &=&\bigvee\limits_{a\in F_1,b\in F_2}f(a)\wedge f(b)
\end{eqnarray*}
Thus,\ $P_L(f)$\ is a\ $\vee-\wedge-$\ morphism.\ $\Box$
\vskip 3mm
Now,\ we will show the main theorem of this section.
\begin{theorem}\
	\begin{enumerate}[(1)]
		\item For each directed space\ $X$,\ each\ $\vee-\wedge-$\ directed space\ $Y$,\ and every morphism\ $f:X\rightarrow Y$,\ there is a unique\ $\vee-\wedge-$\ morphism\ $F$\ from\ $P_L(P_U(X))$\ to\ $Y$\ such that\ $F\circ(\eta_\vee\circ\eta_\wedge)=f$.\ $F$\ is given by\ $P_L(P_U(f))$.\ Thus,\ $P_L(P_U(X))$\ is the free\ $\vee-\wedge-$\ directed space over the directed space\ $X$.
		\item For each directed space\ $X$,\ each\ $\vee-\wedge-$\ directed space\ $Y$,\ and every morphism\ $f:X\rightarrow Y$,\ there is a unique\ $\vee-\wedge-$\ morphism\ $F$\ from\ $P_U(P_L(X))$\ to\ $Y$\ such that\ $F\circ(\eta_\wedge\circ\eta_\vee)=f$.\ $F$\ is given by\ $P_U(P_L(f))$.\ Thus,\ $P_U(P_L(X))$\ is the free\ $\vee-\wedge-$\ directed space over the directed space\ $X$.
		\item For each directed space\ $X$,\ $P_L(P_U(X))$\ and\ $P_U(P_L(X))$\ are isomorphic by a\ $\vee-\wedge-$\ morphism,\ which maps\ $\uparrow(\da x)$\ to\ $\da(\ua x)$.
	\end{enumerate}
\end{theorem}
{\bf Proof}\ \ We prove\ (1)\ and (2)\ is analogous.\ Suppose\ $X$\ is a directed space,\ then\ $P_U(X)$\ is a\ $\wedge-$\ direcrted space and\ $P_U(f)$\ is a\ $\wedge-$\ morphism from\ $P_U(X)$\ to\ $Y$.\ By Theorem\ \ref{comu 1}\ and\ \ref{comu 2},\ $P_L(P_U(X))$\ is a\ $\vee-\wedge-$\ directed space and\ $P_L(P_U(f))$\ is a\ $\vee-\wedge-$\ morphism.\ $f=P_U(f)\circ\eta_\wedge=(P_L(P_U(f))\circ\eta_\vee)\circ\eta_\wedge=F\circ(\eta_\vee\circ\eta_\wedge)$.\ We remain to prove that\ $F$\ is unique.\ Let\ $F_1,F_2$\ be two such\ $\vee-\wedge-$\ morphisms,\ i.e.\ $F_i\circ\eta_\vee\circ\eta_\wedge =f,i=1,2$\ holds.\ We have to show\ $F_1=F_2$.\ Let\ $F_{i}^{'}=F_i\circ\eta_\vee,i=1,2$.\ Since\ $F_{i}^{'}\circ\eta_\wedge=f$,\ by the uniqueness of\ $P_U(f)$,\ $F_1\circ\eta_\vee=F_2\circ\eta_\vee$.\ Since the extension of\ $P_L(P_U(X))$\ is unique,\ we have\ $F_1=F_2$.\ Thus,\ $F$\ is unique,\ and\ $P_L(P_U(X))$\ is the free\ $\vee-\wedge-$\ directed space over\ $X$.\\
(3)\ \ an immediate conclusion by\ (1)\ and\ (2).\ $\Box$

\section*{data availability statement}
Data sharing not applicable to this article as no datasets were generated or analysed during the current study.

\section*{Conflict of interest}
The authors declare that they have no conflict of interest.

\end{document}